\documentclass[preprint,3pt]{elsarticle}

\usepackage{graphicx,amsmath,amsfonts}
\usepackage{times,amsfonts,amssymb,amsmath}
\usepackage{psfrag, color}
\usepackage{multirow}
\usepackage{comment}
\usepackage{caption}
\usepackage{soul}
\usepackage{color}

\usepackage{subfigure}

\usepackage{tikz}
\usetikzlibrary{arrows}

\def\NN{\hbox{I\kern-.2em\hbox{N}}}
\def\RR{{\mathop{{\rm I}\kern-.2em{\rm R}}\nolimits}}
\def\CC{\mathbb{C}}

\def\Q{{\bf Q}}
\def\mg{{\mathcal{G}}}

\def\d{{\bf d}}
\def\f{{\bf F}}
\def\g{{\bf g}}
\def\n{{\bf n}}

\def\r{{\bf r}}
\def\s{{\bf s}}
\def\t{{\bf t}}
\def\u{{\bf u}}
\def\vv{{\bf v}}
\def\x{{\bf x}}
\def\y{{\bf y}}
\def\w{{\bf w}}

\def\z{{\bf z}}

\def\p{{\bf p}}
\def\F{{\bf F}}
\def\f{{\bf f}}
\def\g{{\bf g}}
\def\n{{\bf n}}
\def\x{{\bf x}}
\def\y{{\bf y}}
\def\vv{{\bf v}}
\def\w{{\bf w}}
\def\i{{\bf i}}
\def\j{{\bf j}}
\def\k{{\bf k}}

\def\bp{{\boldsymbol \pi}}

\def\bfax{\mbox{\boldmath$\alpha$}}
\def\bfbx{\mbox{\boldmath$\beta$}}

\def\bnu{\mbox{\boldmath$\nu$}}

\def\bfl{\mbox{\boldmath$\Lambda$}}

\definecolor{orange}{rgb}{1,0.5,0}

\newcommand{\vertiii}[1]{{\left\vert\kern-0.25ex\left\vert\kern-0.25ex\left\vert #1
    \right\vert\kern-0.25ex\right\vert\kern-0.25ex\right\vert}}

\newcommand{\be}{\begin{equation}}
\newcommand{\ee}{\end{equation}}
\newcommand{\ba}{\begin{eqnarray}}
\newcommand{\ea}{\end{eqnarray}}
\newcommand{\supp}{\mathop{\mathrm{supp}}}
\newcommand{\pf}{ \noindent {\bf Proof. }}
\newcommand{\eop} { { \vrule height7pt width7pt depth0pt} \par \bigskip}

\newtheorem{crl}{Corollary}
\newtheorem{rmk}{Remark}
\newtheorem{prn}{Proposition}

\begin{document}

\begin{frontmatter}

\title{ IgA-BEM for 3D Helmholtz problems on multi-patch domains using B-spline tailored numerical integration}
\author[1]{Antonella Falini}

\author[2]{Tadej Kandu\v{c}*}
\cortext[cor1]{Corresponding Author}
\ead{tadej.kanduc@fmf.uni-lj.si}

\author[3]{Maria Lucia Sampoli}

\author[4]{Alessandra Sestini}


\address[1]{Dipartimento di Informatica, Universit\`a degli Studi di Bari Aldo Moro,\\
Via Orabona 4, Bari, Italy}

\address[2]{Faculty of Mathematics and Physics, University of Ljubljana,\\
Jadranska ulica 19,1000 Ljubljana, Slovenia}


\address[3]{Department of Information Engineering and Mathematics, University of Siena,\\
Via Roma 56, Siena, Italy}

\address[4]{Department of Mathematics and Computer Science, University of Firenze,\\
Viale Morgagni 67, Firenze, Italy}

\begin{abstract}
An Isogeometric Boundary Element Method (IgA-BEM) is considered for the numerical solution of Helmholtz problems  on 3D bounded or unbounded domains, admitting a smooth conformal multi-patch representation of their finite boundary surface. The discretization space is formed by $C^0$ inter-patch continuous basis functions whose restriction to a patch simplifies to the span of tensor product B-splines composed with the given patch parameterization. For both regular and singular integration, the proposed model utilizes a numerical procedure defined on the support of each trial B-spline function, which makes possible a function--by--function implementation of the matrix assembly phase. Spline quasi-interpolation is the common ingredient of all the considered quadrature rules; in the singular case  it is combined with  
a B-spline recursion over the spline degree and with a singularity extraction technique, extended to the multi-patch setting for the first time. A threshold selection strategy is proposed to automatically distinguish between nearly singular and regular integrals. Numerical examples on relevant benchmarks show that  the expected convergence orders are achieved with uniform discretization and a small number of uniformly spaced quadrature nodes.
\end{abstract}


\end{frontmatter}
\section{Introduction}\label{sec:Intro}

Isogeometric Analysis (IgA) \cite{LibroHughes} is a powerful tool to obtain a numerical solution of problems governed by partial differential equations, introduced in the literature at the beginning of the new millennium. This new paradigm was motivated by the observation that in engineering applications a domain is described by its boundary parametric representation  generated by Computer Aided Design (CAD) software. For this aim CAD relies on flexible forms, often based on multi-patch tensor product B-spline or rational B-spline (NURBS) spaces, suited also for   applications \cite{Handbook02}. The IgA idea consists in adopting the CAD functional spaces also for approximating the solution of the differential problem, taking into account the well-known optimal approximation power of spline spaces \cite{deBoor01,Schumaker07}. The increasing success of the IgA paradigm in the context of both domain and boundary element methods is easily explained considering that, besides being capable to keep an exact representation of complex domains described in a multi-patch CAD form, IgA formulations guarantee a certain level of accuracy  with considerably less degrees of freedom than traditional Finite Element Analysis (FEA), which relies on larger piecewise polynomial spaces with lower  inter-element smoothness within each patch \cite{Sande2020}. The attractiveness of IgA is even higher for Boundary Element Methods (BEMs), see for example \cite{politis09,Simpson2012,ADSS1,IGABEMbook20,ZANG20201767, Dirk2} and references therein, because, in contrast to domain methods like FEM, they do not require any  preliminary  volumetric parameterization, which can be a remarkably time consuming task. BEMs rely on a boundary integral formulation of the problem, which can be derived whenever the fundamental solution of the associated differential operator is known. Hence, they require only the definition of a mesh  on the boundary of the domain, which is a much easier task especially in the IgA context, where the CAD representation of the boundary is available.

The precursor of this work, dealing with 3D IgA-BEM, is presented in \cite{roma}, where exterior Laplace problems with Dirichlet boundary conditions on single patch domains are considered and, in order to deal also with screen problems, an indirect boundary integral formulation is adopted. In this paper the IgA multi-patch formulation of BEMs (multi-patch IgA-BEMs) based on a direct boundary integral formulation is considered for the numerical solution of problems governed by the Helmholtz equation on 3D bounded or unbounded domains, equipped with either Dirichlet or Neumann boundary conditions. Such problems are of interest in acoustics, since they model radiation and rigid scattering  in the frequency domain \cite{BEMacoustic19}, where radiation refers to the pressure field produced by an object  vibrating within a fluid (usually air or water), while rigid scattering relates to the disturbance caused by an obstacle immersed within an existing acoustic field. The adopted Boundary Integral Equation (BIE) discretized with a collocation approach is the so-called conventional BIE (CBIE), which involves only weakly singular integrals and is augmented with the Sommerfield radiation condition at infinity for exterior problems.  The Helmholtz operator can thus be simply written as a sum of the Laplacian $\Delta$ with an identity operator scaled with a positive parameter $\kappa^2$. However, even when dealing with exterior Helmholtz problems admitting a unique solution, such BIE generates spurious solutions when  $\kappa^2$ is an eigenvalue of $-\Delta$ for the considered boundary conditions. This is the reason why other papers also relying on IgA-BEMs for the numerical solution of the Helmholtz equation introduce its alternative Burton-Miller (BM) boundary integral formulation, which does not suffer from such drawback. Obtained as a linear combination of the mentioned CBIE with the hypersingular BIE (HBIE) also paired  with the Helmholtz operator, the BM integral formulation is clearly more difficult, see for example \cite{Simpson14}, where it is combined with a regularization technique to avoid to deal with hypersingular integrals in the assembly phase. However, successive researches have given numerical evidence that, when $\kappa^2$ is not an eigenvalue of $-\Delta,$ better numerical results are achieved by using the simpler CBIE, see for example \cite{Venas20,Anitescu21}. Since clearly only a  finite number of different values of the parameter $\kappa^2$ can be considered and this can be anyway sufficient for a reasonable reconstruction of an associated time-dependent acoustic field, we have opted for applying the IgA-BEM approach just to the CBIE associated with Helmholtz, with the implicit assumption to avoid values of $\kappa^2$ equal to eigenvalues of $-\Delta.$

The IgA-BEM implementation in this study relies on the B-spline tailored cubature formulas for regular and weakly singular integrals,  introduced in \cite{nash} in the setting of weakly singular surface integrals. The theoretical convergence order of such rules is analyzed in this work for the first time. Since the proposed quadrature schemes are constructed on the support of every bivariate B-spline, a  function-by-function implementation of the matrix assembly phase can be adopted, similarly to work for 2D Laplace problems in \cite{ACDS3}. 
This attractive feature markedly reduces the amount of data accessing in the assembly phase, compared to a more common element--by--element strategy. Another remarkable feature of the proposed quadratures is the freedom of choice for the quadrature nodes -- no reparameterization of nodes near the singularities is needed and, by choosing uniformly spaced nodes in the parametric space for all the required regular and singular numerical integrations, it ensures the possibility of using a small overall number of unique nodes.  
By applying a singularity subtraction technique on singular integrals, the singular and regular part of the integrand are effectively separated. Thanks to the thorough study developed in \cite{tadej21} for the weakly singular kernel expansions, the continuity of the regularized integrands can be controlled and the developed integration rules have consequently better approximation properties for general smooth geometry parameterizations. 
A novelty in this work is a careful decoupling of integrals appearing in the Helmholtz BIE into simpler, real and imaginary, regular and singular integrals. An original strategy for the automatic detection of near singularity is also introduced. Another new key feature lies in coupling the integration routines with multi-patch geometry representations; when the source point and the integration domain are separated by a patch interface, the unmatching parameterization across it poses an additional challenge. To overcome difficulties of this type, for nearly singular inter-patch  integrals a new approximated singularity subtraction technique is introduced. It replaces the source point with its projection on an extrapolation of the patch surface, where the integration has to be performed. 

The paper is organized as follows. The next section introduces the Helmholtz problem and its conventional boundary integral formulation for both  Dirichlet and Neumann boundary conditions. Then, in Section \ref{sec:IgA} the multi-patch IgA setting for the standard CAD representation of a bounded or unbounded domain is described, as well as the finite dimensional functional spaces adopted for the analysis. Section \ref{sec:Quad} describes the developed quadrature formulas, based on the tensor-product formulation of a discrete spline quasi-interpolation approach. The formulas are combined 
with a subtraction regularization technique and applied to the single and double layer Helmholtz kernels. In Section \ref{sec:NE}  the results obtained for some benchmark problems are reported and commented, while some conclusive remarks are given in Section \ref{sec:Conc}. Finally, some technical details related to considered singular kernels are available in Appendix.

\section{The  Helmholtz problem and its boundary integral formulation}\label{sec:Helm}

 We study   3D potential problems described by the  Helmholtz equation with Dirichlet or Neumann boundary conditions on domains $\Omega \subset \RR^3$ admitting a connected closed boundary surface $\Gamma = \partial \Omega,$
\begin{align} \label{BVP}
\left\{ \begin{array}{ll}
\Delta u + \kappa^2 u = 0& {\rm in}\; \Omega,\\
u=u_{\rm D} \quad {\rm or}\; \quad  \frac{\partial u}{\partial \n}= u_{\rm N} & {\rm on}\; \Gamma,
\end{array} \right.
\end{align}
where $u: \Omega \rightarrow \CC$ denotes the unknown {\it potential}, $\n$  the unit normal on $\Gamma$ pointing outward from $\Omega$ and $\kappa>0.$ When the domain is a finite volume the problem is {\it interior} and we use the notation $\Omega^{[i]}$ to underline the type of this domain when it is useful. In the opposite case, we set $\Omega^{[e]}$ and we study the {\it exterior} problem.
Concerning the regularity of the boundary surface $\Gamma,$ for simplicity we assume to be smooth, at least without self-intersections and with a tangent plane well defined at each point and varying continuously.
Note that for exterior problems, equation (\ref{BVP}) has to be augmented with an additional condition, the so called {\it Sommerfeld radiation condition},
\begin{align}
\label{Sommer} \frac{\partial u}{\partial r} - i \kappa u = o \left(\frac{1}{r} \right),
\end{align}
where $r$ denotes the point distance from the origin of the reference system and $i$ is the imaginary unit; see for example \cite{BEMacoustic19}. This condition at infinity is necessary to ensure for any positive $\kappa$ the existence and uniqueness of the solution $u$ for exterior Helmholtz problems for both Dirichlet and Neumann boundary conditions. Uniqueness of $u$ for interior problems is ensured when $\kappa^2$ is not an eigenvalue of the reversed Laplacian operator $-\Delta$ in $\Omega$ \footnote{The scalar $\lambda$ is an eigenvalue for $-\Delta$ if there exists a non vanishing function $u_\lambda$ such that $-\Delta u_\lambda = \lambda u_\lambda$, also fulfilling homogeneous boundary condition of the assigned type. It is possible to prove that the operator $-\Delta$ admits only positive eigenvalues, which define an unbounded infinite sequence of positive numbers depending on the considered  kind of boundary conditions  and on the shape of the finite domain  $\Omega,$ see for example \cite{Laplaceeigen13}.}.
Regarding the regularity of the weak solution $u$ of (\ref{BVP}), we observe that  $u$ belongs to the Sobolev space $H^1(\Omega),$ provided that, for Dirichlet (Neumann) boundary conditions, $u_{\rm D} \in H^{1/2}(\Gamma)\,,\, \left(u_{\rm N} \in H^{-1/2}(\Gamma)\right)$.\\
The Helmholtz equation is of particular interest in acoustic because the solution $u = u(\x,\kappa)$, $\x \in \Omega$ of (\ref{BVP}) can be interpreted as the inverse Fourier transform \cite{Venas19} of the time-dependent scattered pressure field $p= p(\x\,,\,t)$ generated by the scatterer in a given fluid domain, which fulfills the wave equation,
\begin{equation} \label{wave}
\frac{1}{c_f^2} \ \frac{\partial^2 p}{\partial t^2} = \Delta p.
\end{equation}
Namely, setting
$$\breve{p} =\breve{p}(\x\,,\,\nu) := \int_{-\infty}^{+\infty} p(\x\,,\,t) e^{i\nu t}\, dt$$
it is easy to check that $\breve{p}$ verifies the Helmholtz equation with $\kappa = \nu/c_f$ and that, conversely, it is
$$p =\frac{1}{2\pi} \int_{-\infty}^{+\infty} \breve{p}(\x\,,\,\nu) e^{i\nu t}\, d\nu\,.$$
This is of interest because it allows us to deal  with a differential problem formulated just in terms of spatial variables.

In order to approximate the solution of \eqref{BVP} with IgA-BEMs,  we consider the  Conventional Boundary Integral Equation (CBIE) associated with the Helmholtz equation and derived by direct approach,
\begin{align}\label{eqn:BIEdirect}
 \int_{\Gamma} {\mg}_\kappa(\x, \y) \frac{\partial u}{\partial \n} (\y) d\Gamma_{\y} = c(\x) u(\x) + \int_{\Gamma} \frac{\partial {\mg}_\kappa}{\partial \n_{\y}} (\x, \y)\ u(\y) \ d\Gamma_{\y}, \qquad \x \in \Gamma,
\end{align}
where $r : =  \| \r \|_2$ with $\r := \x-\y$. The kernel ${\mg}_\kappa$ and its derivative with respect to the exterior unit normal $\n$ to  $\Gamma$  computed in the point $\y \in \Gamma$ are defined as 
\begin{align*}
{\mg}_\kappa(\x, \y) :=  \frac{1}{4\pi r} e^{{\rm i} \kappa r},  \qquad \qquad  
\frac{\partial {\mg}_\kappa}{\partial \n_{\y}} (\x, \y) \,=\,  \frac{1}{4\pi r} e^{{\rm i} \kappa r} (-\frac{1}{r} + {\rm i} \kappa)\  \frac{\partial r}{\partial \n_{\y}}\,,
\end{align*}
with
\begin{align*}
 \frac{\partial r}{\partial \n_{\y}} = -\frac{\r \cdot \n_{\y}}{r}, \qquad \qquad c(\x) := - \int_\Gamma \frac{\partial \mg_0} {\partial \n_{\y}} (\x,\y) d\Gamma_{\y}\,,
\end{align*}
being $c(\x) \equiv 1/2$ for the assumed smooth surface $\Gamma$. The restrictions of $u$ and $({\partial u}/{\partial \n} ) (\y)$ to $\Gamma$ are usually called {\it Cauchy data}, being one of them coincident with the available boundary datum and the other the unknown.
By separating the real and imaginary parts of $\mg_\kappa,$ we obtain
\begin{align} \label{Gdef}
\mg_\kappa = \frac{1}{4 \pi} \left(\frac{\cos(\kappa r)}{r}+ {\rm i} \frac{\sin(\kappa r)}{r} \right).
\end{align}
Note that the real part of $\mg_\kappa$ goes to infinity as $1/r$ when $r \rightarrow 0^+\,,$ while its imaginary part tends to ${\kappa}/{(4\pi)}$. The other kernel involved in (\ref{eqn:BIEdirect}) is ${\partial {\mg}_\kappa}/{\partial \n_\y}$, which can be written as follows, again separating its real and imaginary parts,
\begin{align} \label{Kdef}
\frac{\partial {\mg}_\kappa}{\partial \n_\y} (\x, \y) = \frac{1}{4\pi}\ \left[ \left( \frac{\r \cdot \n_\y}{r^3} \left(\cos(\kappa r) + \kappa r \sin(\kappa r) \right) \right)\,+\, {\rm i} \frac{\r \cdot \n_\y}{r^2} \left(\frac{\sin(\kappa r)}{r} - \kappa \cos(\kappa r) \right) \right].
\end{align}
For a regular surface $\Gamma$ described by a regular smooth parameterization, the quantity ${\r \cdot \n_\y}/{r^2}$ is bounded for $r \rightarrow 0^+$ but in general it is not continuous at $r = 0$; see the Appendix for the proof. The bounded behavior at $r=0$ of such quantity implies that, analogously to ${\mg}_\kappa$,  the  kernel ${\partial {\mg}_\kappa}/{\partial \n_\y}$ is just weakly singular, see Section 9.1 in \cite{Atkinson97} and also the analytical related insights explicitly reported in the Appendix. Indeed, its real part goes to infinity as $1/r$ when $r \rightarrow 0^+$ while its imaginary part tends to zero because ${\r \cdot \n_\y}/{r^2}$ remains bounded and it is multiplied by a factor going to zero.

Clearly, when the boundary conditions in (\ref{BVP})  are of Dirichlet type, the unknown in (\ref{eqn:BIEdirect})  is the flux $\phi := ({\partial u}/{\partial \n})|\Gamma,$  belonging to  $H^{-1/2}(\Gamma)$  and the CBIE becomes

\begin{align}\label{eqn:BIEdirectD}
 \int_{\Gamma} {\mg}_\kappa (\x, \y)\phi (\y) d\Gamma_{\y} \,=\, c(\x) u_{\rm D}(\x) + \int_{\Gamma} \frac{\partial {\mg}_\kappa}{\partial \n_\y} (\x, \y) u_{\rm D}(\y) d\Gamma_{\y}, \qquad \x \in \Gamma.
\end{align}
Note that this is an integral equation of the general type
\be \label{Symmeq}
V_\kappa \phi(\x) = f(\x),\qquad \x \in \Gamma,
\ee
where $f$ denotes a known function and the single layer operator  $V_\kappa: H^{-1/2}(\Gamma) \rightarrow H^{1/2}(\Gamma)$ is an elliptic isomorphism defined as,
\begin{align*}
V_\kappa \phi(\x) := \int_{\Gamma} \mg_\kappa (\x,\y) \phi(\y)\, d \Gamma_{\y}.
\end{align*}
Conversely, when Neumann conditions are dealt with, the unknown in (\ref{eqn:BIEdirect})  is the acoustic potential $\phi :=u|\Gamma,$  with $\phi \in H^{1/2}(\Gamma)$ and the CBIE becomes
\begin{align}\label{eqn:BIEdirectN}
\int_{\Gamma} \frac{\partial {\mg}_\kappa}{\partial \n_\y} (\x, \y) \phi (\y) \, d\Gamma_{\y} + c(\x) \phi (\x) \,=\,  + \int_{\Gamma} {\mg}_\kappa(\x, \y)u_{\rm N} (\y) d\Gamma_{\y}   \qquad \x \in \Gamma.
\end{align}
This is an integral equation of the general type
\be \label{Symmeq2}
(c(x) I + K_\kappa) \phi(\x) = f(\x),\qquad \x \in \Gamma,
\ee associated with the operator $ c I + K_\kappa,$ where $I$ denotes the identity operator and $K_\kappa: H^{1/2}(\Gamma) \rightarrow H^{1/2}(\Gamma)$ is the following double layer operator,
\begin{align*}
K_\kappa \phi(\x) := \int_{\Gamma}  \frac{\partial {\mg}_\kappa}{\partial \n_\y}(\x,\y)\ \phi(\y)\, d\Gamma_{\y}.
\end{align*}
Note that in our case the right hand side $f$ in both \eqref{Symmeq} and \eqref{Symmeq2} is just the right hand side respectively of \eqref{eqn:BIEdirectD} and \eqref{eqn:BIEdirectN} and also that clearly it is in any case $\phi = \phi(\x, \kappa)$.

If $\phi$ is available, the solution $u = u(\x, \kappa)$, $\x \in \Omega$ of the BVP in (\ref{BVP}) is obtained    by using the so-called {\it representation formula},
\begin{align} \label{rep}
u(\x, \kappa) =  \pm \left( \int_{\Gamma} {\mg}_\kappa(\x, \y)\ \frac{\partial u}{\partial \n} (\y) \ d\Gamma_{\y} - \int_{\Gamma} \frac{\partial {\mg}_\kappa}{\partial \n_\y} (\x, \y) \ u(\y)\ d\Gamma_{\y} \right) ,\qquad \x \in \Omega,
\end{align}
where the sign is positive if $\Omega = \Omega^{[i]}$ (interior problem) and negative otherwise. We note that the numerical implementation of this formula is not completely trivial because near singular integrals have to be approximated, when $\x$ is very near to $\Gamma$. However, one can be merely interested in $\phi$ or in evaluating $u$ just at points sufficiently  far from $\Gamma.$ For example, when an exterior problem is taken into account,  often the interest is just in the {\it far field pattern} of $u,$ that is into the recovery of $u_\infty\,,$ where
\begin{align*}
u_\infty (\w) := \lim_{r \rightarrow \infty} r e^{-i \kappa r} u(r \w)\,, \qquad \w \in S^2,
\end{align*}
where $S^2$ is a unit sphere in $\RR^3$ centered  at the origin of the considered reference system. In such case we can rely on the following formula whose computation just requires rules for regular integrals \cite{Venas20},
\begin{align*}
u_\infty (\w) = \frac{1}{4\pi}\, \int_\Gamma  \left( i \kappa u(\y) (\w \cdot \n(\y)) \,+\,  \frac{\partial u}{\partial \n} (\y)\right) e^{-i \kappa (\w \cdot \y)} d\Gamma_{\y}.
\end{align*}

\section{Multi-patch Isogeometric model}\label{sec:IgA}


In this section we introduce the  IgA-BEM model,  with respect to the geometry representation and the discretization space to express the numerical solution. The flexibility to model complex geometries in the 3D case is facilitated by adopting its multi-patch formulation.
We assume that the boundary ${\Gamma}$ is a union of $M$ patches $\Gamma^{(\ell)}, \ell=1,\ldots,M,$ and for $ \ell \ne k$ it holds $\Gamma^{(\ell)} \cap \Gamma^{(k)} = \emptyset$ and $\partial \Gamma^{(\ell)} \cap \partial \Gamma^{(k)}$ is a boundary edge curve of each patch, a corner point of each
patch or an empty set. To each patch $\Gamma^{(\ell)}$ we assign a geometry map $\F^{(\ell)}$, such that   $\overline{\Gamma}^{(\ell)} = \mbox{Image}(\F^{(\ell)})$,  $\F^{(\ell)}~:~ [0,1]^2  \rightarrow \overline\Gamma^{(\ell)} \subset \RR^3$ is regular and belongs to $C^2([0,1]^2)$. 
Each $\F^{(\ell)}$ is written in the following standard NURBS form \cite{Handbook02},
\begin{align} \label{NURBS} \F^{(\ell)}(\t) := \frac{ \displaystyle{\sum_{\i \in {\cal J_{\rm g}}^{(\ell)} }} w_{\i}^{(\ell)}\, \Q_{\i}^{(\ell)}\, \hat B_{\i,\d_{\rm g}}^{(\ell)}(\t) }{\displaystyle{\sum_{\i \in {\cal J_{\rm g}}^{(\ell)}} } w_{\i}^{(\ell)} \, \hat B_{\i,\d_{\rm g}}^{(\ell)} (\t)}, \quad \t \in [0,1]^2,
\end{align}
 where a set $\{ \Q_{\i}^{(\ell)} \in \RR^3, \i \in {\cal J_{\rm g}}^{(\ell)} \}$  defines a net of  control points which, together with the associated set of positive weights $\{w_{\i}^{(\ell)} \in \RR,  \i \in {\cal J_{\rm g}}^{(\ell)}\},$ is the basic element typically used in the CAD environment to design free-form surfaces. The tensor product B-spline basis function $\hat B_{\i,\d_{\rm g}}^{(\ell)}(\t) := \hat B_{i_1,d_{\rm g,1}}^{(\ell)}(t_1)\  \hat B_{i_2,d_{\rm g,2}}^{(\ell)}(t_2)$, where $\i:=(i_1,i_2) \in \cal J_{\rm g}^{(\ell)}$, is of bi-degree $\d_{\rm g}=(d_{\rm g,1},d_{\rm g,2})$, defined on  $[0,1]^2$ in variable $\t :=(t_1,t_2)$ and with respect to the bidimensional knot array $T_{\rm g}^{(\ell)}$ associated with the $\ell$-th patch, $T_{\rm g}^{(\ell)} := T_{\rm g, 1}^{(\ell)} \times T_{\rm g,2}^{(\ell)}$, and  $T_{\rm g,1 }^{(\ell)}$ and $T_{{\rm g},2}^{(\ell)}$ are both clamped knot vectors with entries in the interval $[0,1]$. When two patches share an edge ${\cal C}_{\ell,k}:=\partial {\Gamma}^{(\ell)} \cap \partial {\Gamma}^{(k)}$ for $\ell \neq k$, the shared edge is a NURBS curve.  For simplicity we usually assume it to be parameterized with the same geometry map for both patches, up to reversion in directions and direction swapping in the parametric space.

As commonly done also in the literature, in order to gain more flexibility, in our generalized IgA-BEM setting the unknown Cauchy datum $\phi$ is approximated in a functional space whose restriction to $\bar \Gamma^{(\ell)}$ is the geometric counterpart of  the space generated from bivariate B-splines always defined in $[0, 1]^2$ but not necessarily obtained with a refinement procedure from $\hat B_{\i,\d_{\rm g}}^{(\ell)},\, \i \in {\cal J_{\rm g}}^{(\ell)} $. Such space is characterized by its own degree $\d = (d_1,d_2)$ (common to all patches) and its extended knot vectors, ${T}_j^{(\ell)}, j=1,2.$ 
Then, denoting with   $ h_j^{(\ell)}$  the maximal distance between successive knots in ${T}_j^{(\ell)},$ setting $h^{(\ell)}:= \max\{ h_1^{(\ell)},\, h_2^{(\ell)}\}$ and introducing a set of bivariate indices ${\cal J}^{(\ell)}$ for the basis functions, we can introduce $\hat {\cal S}^{(\ell)}_{h^{(\ell)}} $ which denotes the  spline space generated by the tensor product B-spline basis $\{ \hat B_{\j,\d}^{(\ell)}:  \j = (j_1,j_2) \in {\cal J}^{(\ell)} \}$ of bi-degree $\d$ defined on $[0,1]^2$ with respect to the bidimensional knot array ${T}^{(\ell)} := {T}_1^{(\ell)} \times { T}_2^{(\ell)}$.
\subsection{Discontinuous basis}
When no continuity constraint is imposed across patches for the basis functions, the global   space used to approximate the unknown Cauchy datum $\phi$ is
\begin{align*}
 {\cal S}_{h,\d} \,:=\, \mathrm{span}\{ { B}_{\j,\d}^{(\ell)}: \j \in {\cal J}^{(\ell)}, \quad  \ell=1,\ldots,M\},
\end{align*}
where $h := \max\{ h^{(\ell)}, \ell=1,\ldots,M\}$ and
\begin{align*}
{B}_{\j,\d}^{(\ell)}(\x) =  \left\{
\begin{array}{cc}
\hat B_{\j,\d}^{(\ell)} \circ {\F^{(\ell)}}^{-1}(\x), & \mbox{if }  \x \in {\Gamma}^{(\ell)}, \cr
 0, & \mbox{otherwise.}
 \end{array}
 \right.
\end{align*}
 Clearly with this setting the global number $N_{\rm DOF}$ of degrees of freedom used to approximate $\phi$ is 
\begin{align*}
{N}_{\rm DOF} = \sum_{\ell =1}^M \vert {\cal J}^{(\ell)} \vert,
\end{align*}
where $ \vert {\cal J}^{(\ell)} \vert$ is the cardinality of the set.
 In this setting it is convenient to select $ \vert \cal J^{(\ell)} \vert$ distinct collocation points belonging to the interior of each patch $\Gamma^{(\ell)}$ , $\ell=1,\ldots,M.$ In particular we fix $\x_\j^{(\ell)} = \F^{(\ell)} (\s_\j)$, $\j \in {\cal J}^{(\ell)},$ where $\{ \s_\j^{(\ell)}: \j \in {\cal J}_\ell \}$ denotes a set of points in $(0,1)^2$, defined as the Cartesian product of two sets of abscissas, the so-called improved Greville abscissas associated with ${T}_1^{(\ell)}$ and ${T}_2^{(\ell)}$ \cite{WangBenson15}. We recall that the improved Greville abscissae related to a set of Greville abscissae $\bar \xi_1,\bar \xi_2, \ldots,\bar \xi_n$, are Greville abscissae whose first and  last elements are modified as follows:
\begin{equation} \label{imp}
 \bar \xi_1=\bar \xi_1+\omega(\bar \xi_2-\bar \xi_1) ,\;\;\; \bar \xi_n=\bar \xi_n-\omega(\bar \xi_n-\bar \xi_{n-1})\,, \quad \mbox{with} \,\, 0<\omega<1\,.
\end{equation}

Thus, a discrete version of the integral formulation of the Dirichlet \eqref{eqn:BIEdirectD} or Neumann \eqref{eqn:BIEdirectN} problem is obtained by approximating
$\phi$ in the finite dimensional composite space ${\cal S}_{h,\d}$.
The applied collocation method leads to a linear system
\begin{align} \label{sistlin}
A \,\bfax = \bfbx.
\end{align}
The unknown entries in the vector $\bfax$ are the coefficients which allow us to define our approximation $\phi_h$ of $\phi$ in the space ${\cal S}_{h,\d}$ patch-wisely as follows,
\begin{align*}
\phi_h(\x)  := \sum_{\ell=1}^M \sum_{\j \in {\cal J}^{(\ell)}} \alpha^{(\ell)}_\j {B}_{\j,\d}^{(\ell)}(\x), \qquad  \x \in {\Gamma}.
\end{align*}
{
The square system matrix $A \in \mathbb{C}^{{ N_{\rm DOF}} \times N_{\rm DOF}}$ and the right-hand side vector $\bfbx \in \mathbb{C}^{N_{\rm DOF}}$ have a block representation,  $A = (A^{ (\ell,k) } ),\,\, \bfbx = (\bfbx^{ (\ell) })$, for $\ell, k =1,\dots,M$, where  $A^{(\ell,k)} \in \mathbb{C}^{ \vert {\cal J}_\ell \vert \times \vert {\cal J}_k \vert }$ and $\bfbx^{(\ell)} \in \mathbb{C}^{ \vert {\cal J}_\ell \vert}$.
For each pair of patches $\Gamma^\ell,\Gamma^k$, the rows  (columns) of the related submatrix are ordered by using a lexicographical ordering of the elements in $\mathcal J^{(\ell)}$ ($\mathcal J^{(k)}$), namely we identify a single index $i$ with $\i$ and $j$ with $\j$. After this simplification, for the Dirichlet case such entries read as follows in both the physical and the parametric domain,
\begin{align}
\label{eqn:DL}
A^{(\ell,k)}_{i,j} := &\int_{\Gamma^{(k)}}  {\mg}_{\kappa}(\x_\i^{(\ell)},\y)  B^{(k)}_{\j,\d} (\y)\   d \Gamma
\,= \int_{[0,1]^2}  {\mg}_{\kappa}(\x_\i^{(\ell)}, \F^{(k)} (\t)) \hat B^{(k)}_{\j,\d} (\t)\  J^{(k)} (\t) \ d \t,
  \\
\label{eqn:DR}
\bfbx^{(\ell)}_i := & \sum_{k=1}^M \int_{\Gamma^{(k)}} \frac{\partial \mg_\kappa}{\partial \n_\y}(\x_\i^{(\ell)},\y)\ u_{\rm D} (\y)\ d \Gamma +  \frac{1}{2} u_{\rm D} (\x_\i^{(\ell)}) = \nonumber\\
\ &  \sum_{k=1}^M \ \int_{[0,1]^2} \frac{\partial \mg_\kappa}{\partial \n_\y}(\x_\i^{(\ell)},\F^{(k)}(\t))\ u_{\rm D} (\F^{(k)}(\t))\  J^{(k)}(\t) \ d \t + \frac{1}{2} u_{\rm D} (\x_\i^{(\ell)}),
\end{align}
where $i=1,\ldots, \vert{\cal J}^{(\ell)}\vert, \quad j=1,\ldots, \vert{\cal J}^{(k)}\vert,\,\, \hat B^{(k)}_{\j,\d} \in {\hat {\cal S}}^{(k)}_{h^{(k)}}$ and $J^{(k)}$ represents the infinitesimal surface area  element on the $k$-th patch,
$$J^{(k)}(\cdot) := \left \Vert \frac{\partial \F^{(k)}}{\partial t_1}(\cdot) \times   \frac{\partial \F^{(k)}}{\partial t_2}(\cdot) \right \Vert_2.$$
Conversely, for the Neumann case we have
\begin{align}
\label{eqn:NL}
A^{(\ell,k)}_{i,j} = &:= \int_{\Gamma^{(k)}} \frac{\partial \mg_\kappa}{\partial \n_\y}(\x_\i^{(\ell)},\y)  B_{\j,\d}^{(k)}(\y) \ d \Gamma + \frac{1}{2} B_{\j,\d}^{(k)}(\x_\i^{(\ell)}) =  \nonumber \\
\ & \int_{[0,1]^2} \frac{\partial \mg_\kappa}{\partial \n_\y}(\x_\i^{(\ell)}, \F^{(k)}(\t)) \hat B_{\j,\d}^{(k)}(\t) J^{(k)}(\t)\ d \t + \frac{1}{2} \hat B_{\j,\d}^{(k)}(\x_\i^{(\ell)}),
 \\
\label{eqn:NR}
\bfbx^{(\ell)}_i := & \sum_{k=1}^M  \ \int_{\Gamma^{(k)}} \mg_\kappa(\x_\i^{(\ell)}, \y)\ u_{\rm N} (\y)\  d \Gamma
\,= \sum_{k=1}^M \ \int_{[0,1]^2} \mg_\kappa(\x_\i^{(\ell)}, \F^{(k)}(\t))\ u_{\rm N} (\F^{(k)}(\t))\  J^{(k)}(\t)\ d \t,
\end{align}
Referring for brevity just to the Dirichlet case, we report below an equivalent scaled expression of these entries, since it  will  be useful in Section \ref{sec:Quad} to develop a reasonably scaled  procedure to detect near singular integrals to be computed during the assembly phase. Denoting with $\mu_k$ a reference length for the $k$-th patch (e.g., the square root of its area or its diameter), with $\varphi_k$ the associated uniform scaling such that $\tilde \z = \varphi_k(\z) = \z/ \mu_k$ and with $\tilde{\Gamma}^{(k)}$ the conformally scaled version of the $k$--th patch, we have also
\begin{align}
\label{eqn:DL2}
A^{(\ell,k)}_{i,j} = & \mu_k \ \int_{\tilde{\Gamma}^{(k)}}  \mg_{\kappa_k}(\varphi_k(\x_\i^{(\ell)}), \tilde{\y})  B^{(k)}_{\j,\d} (\mu_k \tilde{\y})\   d \tilde{\Gamma} = \nonumber \\
\ & \mu_k \  \int_{[0,1]^2}  {\mg}_{\kappa_k}(\varphi_k(\x_\i^{(\ell)}), \tilde{\F}^{(k)} (\t)) \hat B^{(k)}_{\j,\d} (\t)\  \tilde{J}^{(k)} (\t) \ d \t,
 \\
\label{eqn:DR2}
\bfbx^{(\ell)}_i = & \sum_{k=1}^M \ \int_{\tilde{\Gamma}^{(k)}} \frac{\partial {\mg}_{\kappa^{(k)}}}{\partial \n_{\tilde y}}(\varphi_k (\x_\i^{(\ell)}), \tilde \y)\ u_{\rm D} (\mu_k \tilde \y)\ d \tilde{\Gamma} + \frac{1}{2} u_{\rm D} (\x_\i^{(\ell)}) = \nonumber \\
\ & \sum_{k=1}^M \ \int_{[0,1]^2} \frac{\partial \mg_{\kappa_k}}{\partial \n_{\tilde{\y}}} (\varphi_k(\x_\i^{(\ell)}),\tilde{\F}^{(k)}(\t))\ u_{\rm D} (\mu_k \tilde{\F}^{(k)}(\t))\  \tilde{J}^{(k)}(\t) \ d \t + \frac{1}{2} u_{\rm D} (\x_\i^{(\ell)}),
\end{align}
where $\tilde \F^{(k)} := \F^{(k)}/\mu_k,\, \tilde{J}^{(k)}$ is the related infinitesimal scaled area element  and $\kappa_k:= \kappa \mu_k.$

\subsection{$C^0$ continuous basis}

Let us denote the space of B-spline functions that join with $C^0$ regularity across patches by ${\cal S}^0_{h,\d}: = {\cal S}_{h,\d} \cap  C^0({\Gamma})$. For simplicity let us assume that any two adjacent patches $\Gamma^{(\ell)}, \Gamma^{(k)}$ with a common boundary curve ${\cal C}_{\ell,k}$ have the same knot vector on the edges of the patches. Namely, on $\ell$-th patch we consider the knot vector $T_1^{(\ell)}$ if $\F^{(\ell)}(t_1,0)$ or $\F^{(\ell)}(t_1,1)$ define the curve ${\cal C}_{\ell,k}$, and $T_2^{(\ell)}$ in the other case. The same procedure is done for $k$-th patch. Our assumption then states that the two vectors $T_i^{(\ell)}$, $T_j^{(k)}$ need to coincide (up to a possible reversion of one knot sequence, e.g., when $\F^{(\ell)}(1,t_2) = \F^{(k)}(1,1-t_2)$). For a more general setting to describe the geometry and the discretization space without the matching of adjacent knot vectors we refer to \cite{LibroHughes}.

The so-called interior basis functions for ${\cal S}^0_{h,\d}$ includes all basis functions ${B}_{\j,\d}^{(\ell)}$ from ${\cal S}_{h,\d}$, whose corresponding B-splines $\hat{B}_{\j,\d}^{(\ell)}$ vanish on $\partial( [0,1]^2)$. The so-called edge or vertex basis functions are obtained by identifying and summing together the remaining basis functions from ${\cal S}_{h,\d}$ with shared knot vectors. More precisely, an edge basis function $B\in {\cal S}^0_{h,\d}$ associated to ${\cal C}_{\ell,k}$ is a sum of two remaining (not interior) functions $B_{\i,\d}^{(\ell)}$ and $B_{\j,\d}^{(k)}$ from ${\cal S}_{h,\d}$, that vanish at the endpoints of ${\cal C}_{\ell,k}$ and satisfy
$\hat B_{\i,\d}^{(\ell)} \circ {\F^{(\ell)}}^{-1}(\x) =  \hat B_{\j,\d}^{(k)} \circ {\F^{(k)}}^{-1}(\x)$ for every $\x \in {\cal C}_{\ell,k}$. Similarly, a vertex function associated to a patch corner vertex $\x_v$, is defined as a sum of all
 basis functions $B_{\i,\d}^{(\ell)}$ from ${\cal S}_{h,\d}$, $\ell=1,\dots,M$ and $\i \in \cal J^{(\ell)}$ , that satisfy $\hat B_{\i,\d}^{(\ell)} \circ {\F^{(\ell)}}^{-1}(\x_v)\neq 0$.

In all our experiments with the $C^0$ continuous basis we have fixed the collocation points as a set of images of Cartesian products of standard Greville abscissas  on $[0,1]^2$, mapped to $\overline \Gamma^{(\ell)}$ via map $\F^{(\ell)}$ for all $\ell=1,\dots,M$, and by counting just once the repeated points, which appear on the boundary curves of patches. 

The approximate solution $\phi_h$ is constructed similarly to the discontinuous case, bearing in mind a variation in the collocation points and a sum of adequate entries in
\eqref{eqn:DL}, \eqref{eqn:NL} over all suitable $k$ to evaluate integrals that span over the patches $\Gamma^{(k)}$ adjacent to $\Gamma^{(\ell)}$.


\section{Quadratures} \label{sec:Quad}

\subsection{Singularity extraction}
Let us  introduce our novel numerical integration method  which makes possible a  function--by--function assembly phase. For the sake of simplicity, we refer  just to the one patch setting, i.e., we assume $M=1,$ thus avoiding any reference to the patch index, when the source point and the integration domain belong to the same patch. Some considerations about integration in the more general multi-patch setting are also given.

It is convenient to write the two singular kernels $\mg_\kappa$ and ${\partial {\mg}_\kappa}/{\partial \n_\y}$ as in \eqref{Gdef} and \eqref{Kdef} and analyze which terms are singular, hence a singular integration has to be considered, and which  terms define regular integrals. Referring to the Appendix for some explicit insights, as already noted in Section \ref{sec:Helm},  the real part in \eqref{Gdef} is weakly singular at $r=0$ with its non-singular part $\cos(\kappa r)$  smooth under the assumption of a smooth geometry parameterization. For the same reason, its imaginary part  is a smooth function. Similarly,  the  factor $\r \cdot \n_\y /r^3$ of the real part in \eqref{Kdef} is weakly singular, while the other factor, $(\cos(\kappa r) + \kappa \sin(\kappa r) r),$ is smooth. Finally,  the imaginary part in \eqref{Kdef} is in general only $C^1$ continuous at $r=0$ because the factor ${\r \cdot \n_\y} / {r^2}$ is bounded and $(\sin(\kappa r) / r - \kappa \cos(\kappa r))$ has a double zero at $r=0$.
In our numerical experiments such smoothness of the integrands was sufficient to numerically integrate the corresponding integrals with the required accuracy.

\begin{rmk}
To improve accuracy of the numerical integration for integrals involving the imaginary part of \eqref{Kdef}, a more sophisticated approach could be considered. For example, the integrals could be handled similarly as the singular ones, treating the part ${\r \cdot \n_\y} / {r^2}$ as the ``less regular part'' of the integrand, since the remaining part $\left({\sin(\kappa r)}/r - \kappa \cos(\kappa r) \right)$ is an analytic function.
\end{rmk}

In order to give a compact introduction of our quadrature approach, first we denote with $\s \in [0, 1]^2$ the pre-image in the parametric domain of a fixed generic collocation point $\x$. Note that for brevity in the following we refer to the patch index only when it is essential. Conversely, we refer to the patch scaled integrals appearing for example in formulas \eqref{eqn:DL2}--\eqref{eqn:DR2}, as well as to the related notation also introduced in Section \ref{sec:IgA}. So for example we denote with $\tilde \x $ the scaling of $\x $, assuming that the considered scaling is  that associated with the patch where the integration is developed. Then,  we can introduce the following general form representing any  integral to be dealt with,
\begin{align} \label{genint}
\int_{R_\j}U (\s, \t) \hat B_{\j,\d} (\t)\ g(\t)\,d\t,
\end{align}
where  $g$ is a given smooth  function, $\hat B_{\j,\d}$ is a B-spline basis function, $R_\j := \supp \{ \hat B_{\j,\d} \}$ and $U$ refers to any of the following two weakly singular kernels respectively associated with the single and double layer potential,
\begin{align}\label{eqn:kernelU}
U(\s , \t) = U_{\rm SL}(\s, \t) :=  \frac{1}{\| \tilde \x(\s) - \tilde \F(\t)\|_2}
\qquad \textrm{and} \qquad
U(\s , \t) = U_{\rm DL}(\s, \t) := {\frac{(\tilde \x (\s) - \tilde\F(\t)) \cdot \tilde \bnu (\t) }{\| \tilde \x(\s) - \tilde \F(\t)\|_2^3}},
\end{align}
with $\tilde \bnu(\t) := ({\partial \tilde \F} / {\partial t_1})(\t) \times ({\partial \tilde \F} / {\partial t_2})(\t) $ and $ \tilde \F$ denoting the scaled normal and the scaled mapping associated to the conformally scaled surface patch. 

\begin{rmk}
Integrals defining the entries of $\bfbx$ do not include the basis functions $\hat B_{\j,\d}$ directly. To split the integration domain $[0,1]^2$ into smaller subdomains and  arrive to the above general form, we introduce  additional knot vectors, $T_{{\rm r},1}, T_{{\rm r},2}$, and inside (\ref{eqn:DR}) and (\ref{eqn:NR}) we add a sum of all B-spline basis functions $\hat B_{\j,\boldsymbol 0}$ of degree zero, defined on $T_{{\rm r},1} \times T_{{\rm r},2}$. Therefore, the approach to evaluate integrals in  $\bfbx$ is analogous to the one for the system matrix $A$.
\end{rmk}
It is easy to determine when the integral in \eqref{genint} is actually singular. Indeed, since we are dealing with smooth surfaces without self-intersections, this occurs if and only if  $\tilde \x$ belongs to the scaled patch where the integration is developed, and the associated $\s$ belongs to $R_\j$.
On the other hand, detection of nearly singular integrals that require analogous treatment is more delicate if we want to preserve the accuracy of the numerical scheme without deteriorating its efficiency. 
In such setting we propose to consider nearly singular the integral in \eqref{genint} when 
\begin{align} \label{eqn:thresh_dist}
\tilde r_{{\rm min},\j}(\x):= \min_{\t  \in R_\j} \| \tilde \x - \tilde \F(\t) \|_2 \leq \,  c \ \delta,
\end{align}
where $0<c<1$ is a user-defined constant and $\delta$ is a suitable patch-depended threshold. 
For further details concerning the definition of $\delta$ see the discussion in Subection~\ref{thresh}.

When the condition \eqref{eqn:thresh_dist} is not satisfied, the kernel $U$ is considered regular and integration rules for regular integrals are applied. Conversely, when the kernel $U$ in (\ref{genint}) is weakly singular or near singular, we preliminarily adopt a regularization technique based on the subtraction of singularity,
\begin{align}\label{eq:subtraction}
\int_{R_\j}U(\s , \t) \hat B_{\j,\d}(\t) g(\t)\,d\t =
\int_{R_\j} \big( U (\s , \t)-{U}_{\s}^m(\s-\t) \big) \hat B_{\j,\d}(\t) g(\t)\,d\t + \int_{R_\j}{U}_{\s}^m(\s-\t) \hat B_{\j,\d} (\t) g(\t)\,d\t,
\end{align}
where ${U}_{\s}^m$ is an approximation of $U(\s,\bullet)$, obtained by truncating a particular series expansion of $U$ about the singular point $\t=\s$ after the $m$-th term; see \cite{nash} for the first application of these expansions in the simplest form and \cite{tadej21} for detailed construction and analysis.
Note that the singularity subtraction technique that splits an integral into two ones is applied only to a small portion of integrals. In this case, the first term of the right hand side of \eqref{eq:subtraction} becomes a \emph{regularized integral} while the second one is a singular integral with a simplified singular kernel. The continuity of the integrand in the regularized integral is controlled by the number $m$ of terms in $U_{\s}^m$ and is equal to $C^{m-2}$ at $\t=\s$ for the chosen $m\geq 1$ (with $C^{-1}$ we denote functions that are integrable in the standard sense but have discontinuities of the first kind at $\t=\s$), if the geometry map $\F$ is a smooth enough function \cite{tadej21}. 

By defining local coordinates $\z:=\s -\t$, the kernel reads
\begin{align}\label{eqn:kernel}
U^m_\s (\z)  = \sum_{k=1}^m R_\s(\z)^{-2 (k+ \gamma)+1}\; P_{\s,3k + 2 \gamma - 3}^{[k]} (\z),
\end{align}
where $R_\s(\z) := \sqrt{\z^\top M_\s \z}$ is a square root of a bivariate quadratic homogeneous
polynomial and $M_\s$ is the $2 \times 2$ symmetric positive definite matrix associated with the first fundamental form of a smooth surface (e.g., see formulas (\ref{EFG}) and (\ref{matM}) in the Appendix), 
$P_{\s,3k + 2\gamma - 3}^{[k]}$ are suitable homogeneous polynomials of degrees $3k + 2 \gamma - 3$, and  $\gamma=0$ and $\gamma=1$ for  $\mg_\kappa$ and  ${\partial {\mg}_\kappa}/{\partial \n_\y}$, respectively. Note that coefficients inside $R_\s$ and $P_{\s, 3k + 2\gamma - 3}^{[k]}$ depend on a chosen $\s$ and on derivatives of the map $\tilde \F$.
For example, for $m=1$ the kernels of the single and double layer potentials are  approximated as follows,
\begin{align*}
U^1_\s (\z)  = \left\{ \begin{array}{ll} \frac{1}{4\pi} \left[ R_\s(\z)^{-1}  \right] & \mbox{if } U = U_{\rm SL} \,,\cr  
 \frac{1}{4\pi} \left[ R_\s(\z)^{-3} P^{[1]}_{\s,2}(\z)  \right] & \mbox{if } U = U_{\rm DL}\,, \cr
 \end{array} \right.  
\end{align*}
where $P^{[1]}_{\s,2} (\z) = \frac{1}{2} \left(L \ z_1^2  + 2M \ z_1 z_2+ N z_2^2 \right)$,
and $L,M,N$ are the coefficients of the second fundamental form of the geometry, evaluated at $\s$, see \eqref{LMN} in the Appendix.
\bigskip

Relating to the more general multi-patch setting, when source point $\x=\F^{(\ell)}(\s)$ and integration variable $\y = \F^{(k)}(\t)$ belong to different adjacent patches, $\Gamma^{(\ell)}$ and $\Gamma^{(k)}$, the integral in \eqref{genint} is nearly singular if condition \eqref{eqn:thresh_dist} is satisfied, and regular otherwise. In the latter case, the procedure is that adopted for regular integrals on one patch while the former case clearly poses an additional challenge. 
In order to explain our procedure in such setting,  we need first to make the notation more precise for $\ell \neq k$, adding the superscript $(\ell,k)$ when referring to the kernel $U,\ U \rightarrow U^{(\ell,k)},$ whenever $U_{\rm SL}$ or $U_{\rm DL}$ is referred to,
\begin{align}\label{eqn:acrossPatchU}
 U^{(\ell,k)}_{\rm SL}(\s, \t) :=  \frac{1}{\| \tilde \F^{(\ell)}(\s) - \tilde \F^{(k)}(\t)\|_2},
\quad 
 U^{(\ell,k)}_{\rm DL}(\s, \t) := {\frac{(\tilde \F^{(\ell)}(\s) - \tilde \F^{(k)}(\t)) \cdot \tilde \bnu^{(k)}(\t)}{\| \tilde \F^{(\ell)}(\s) - \tilde \F^{(k)}(\t)\|_2^3}},
\end{align}
indicating that the source point $\tilde \x$ is connected with $\s$ in the parametric domain via the geometry map $\tilde \F^{(\ell)}$ but integration has to be developed in the $k$-th patch. 
Since the parameterization of the geometry is prescribed patch-wisely, the proposed series expansion $U_\s^{(\ell,k),m}$ is limited to the $\ell$-th patch. On the other hand, since $\F^{(k)}$ is regular, we can (slightly) enlarge its definition outside the original domain $[0,1]^2$ and still have a regular geometry map; for convenience let us use the same notation for the extended function,
$\F^{(k)}: [0,1]^2 \subset D_k  \rightarrow  \RR^3$. Clearly on the enlarged part of the domain the map $\F^{(k)}$ describes just an approximation of the actual geometry, which is exactly described by $\F^{(\ell)}$. By replacing $\tilde \x = \tilde \F^{(\ell)}(\s)$ in \eqref{eqn:acrossPatchU} by $\tilde \F^{(k)}(\s_e)$ with
\begin{align*}
\s_e :=\arg \min_{\s' \in D_k}  \| \tilde \x - \tilde \F^{(k)}(\s')\|_2\,,
\end{align*}
we obtain the kernel $U_{\s_e}^{(k,k),m}$ that approximates $U_{\s}^{(\ell,k),m}$ and naturally generalizes the case $\ell=k$. Since the actual source point $\tilde \x$ is replaced by its nearby point $\tilde \F^{(k)}(\s_e)$, the cancellation of the singularity is not exact. This does not present any real problem because only the evaluation of inter-patch nearly singular integrals is required.

As an example, let $\Gamma^{(\ell)}$ and $\Gamma^{(k)}$ be two adjacent NURBS patches of a unit sphere  (see the beginning of Section~\ref{sec:NE} for more details and Figure~\ref{fig:patch}(a) for visualization of the sphere representation). In Figure~\ref{fig:extrapolation} the function  $\t \mapsto | U^{(\ell,k)}(\s, \t) - U_{\s_e}^{(k,k),2} (\s_e - \t) | $ is depicted for a chosen position of $\x = \F^{(\ell)}(\s) \in \Gamma^{(\ell)}$ near the interface between $\Gamma^{(\ell)}$ and $\Gamma^{(k)}$ (for simplicity, here we do not consider any scaling, that is we assume $\mu_k=1$). Such function is the absolute value of the regularized part of the integrand in \eqref{eq:subtraction} for the kernels of the single and double layer potential. The figure shows its value in the whole parametric domain $[0,1]^2$ of $\Gamma^{(k)}\,,$ together with the position of $\s_e$ in $D_k.$ The distance between the computed nearby point $\F^{(k)}(\s_e)$ and the actual source point $ \x = \F^{(\ell)}(\s)$ is approximately $4.1 \cdot 10^{-9}$ and $9.5 \cdot 10^{-9}$, when the source point is close to a vertex on the interface  (left figures) and to its central part (right figures), respectively. The figure shows that the regularized function $U^{(\ell,k)}(\s, \t) - U_{\s_e}^{(k,k),2} (\s_e - \t)$ has small oscillations and is bounded--in comparison consider that the maximum value of the function $\t \mapsto | U^{(\ell,k)}(\s, \t)|$ in $[0,1]^2$ is $41.9$  and $35.2$ for the two cases with the single layer potential (top figures) and $25.0$, $31.6$ for the cases with the double layer potential (bottom figures). Thus, it can be effectively integrated with the quadrature rule for regular integrals presented in Subsection~\ref{sec:regularRule}.

\begin{figure}[thb]
\centering
{\includegraphics[trim = 4cm 9.25cm 4.5cm 10cm, clip = true, height=5.5cm]{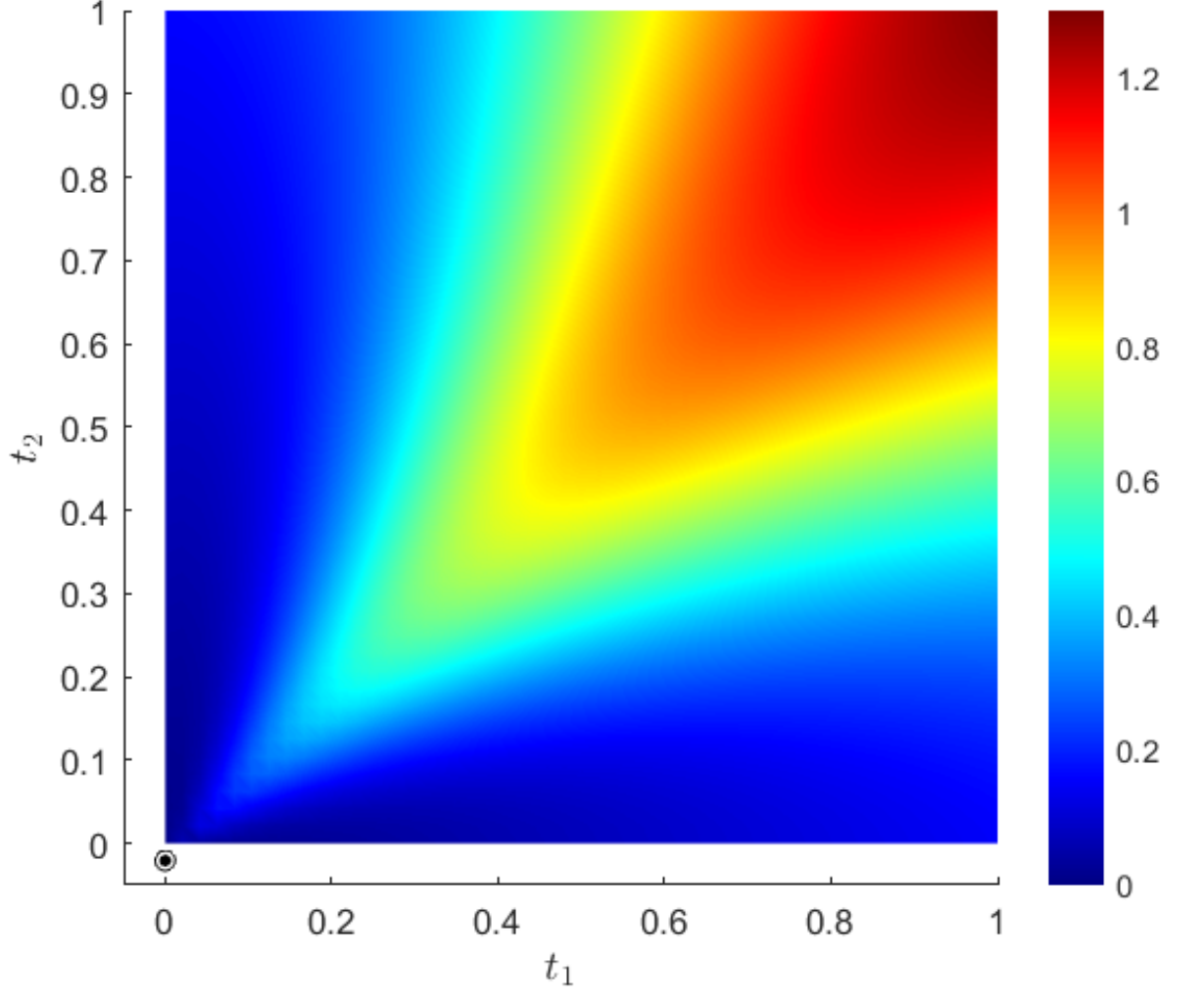}}
{\includegraphics[trim = 4cm 9.25cm 4.5cm 10cm, clip = true, height=5.5cm]{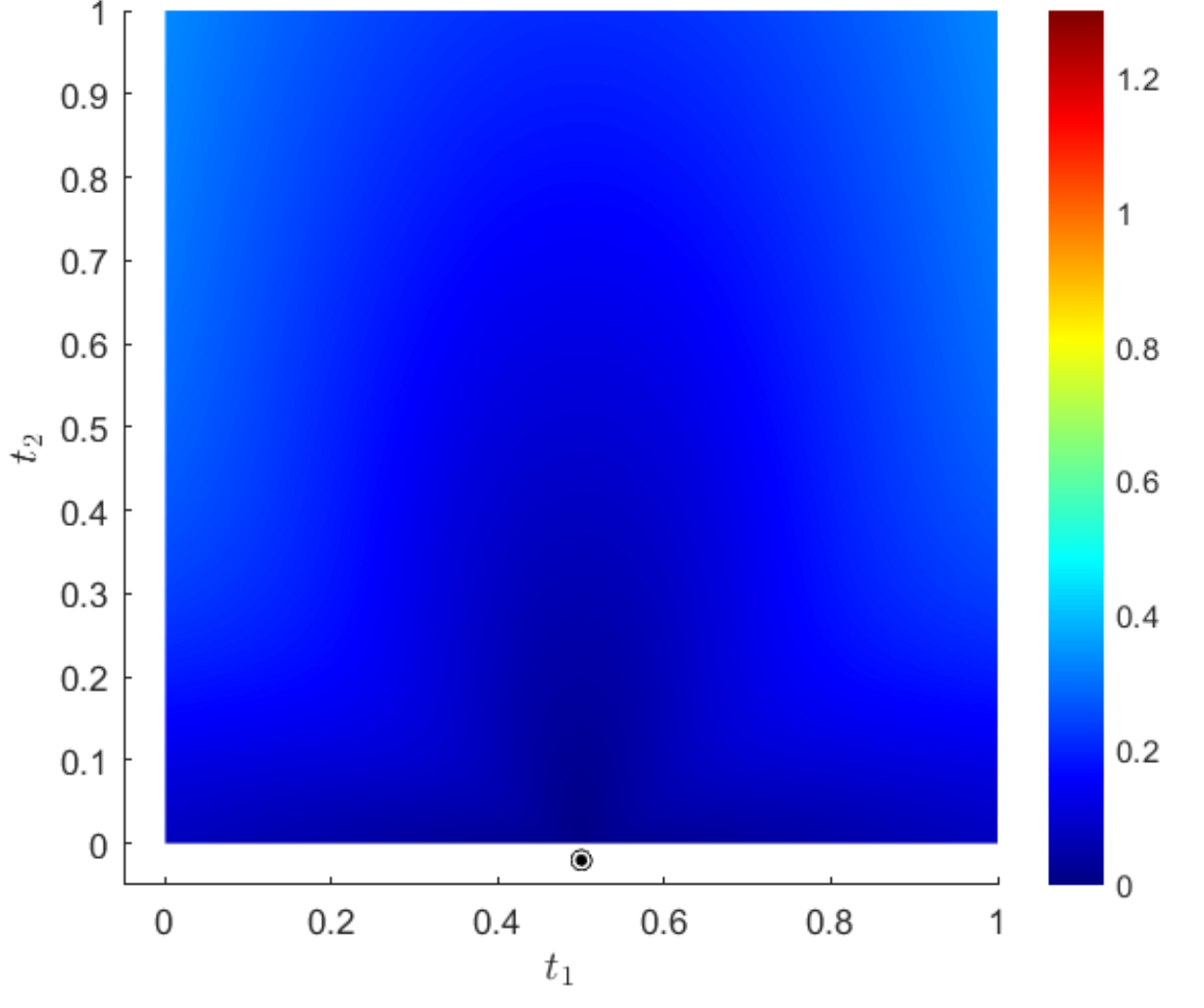}}
{\includegraphics[trim = 4cm 9.25cm 4.5cm 10cm, clip = true, height=5.5cm]{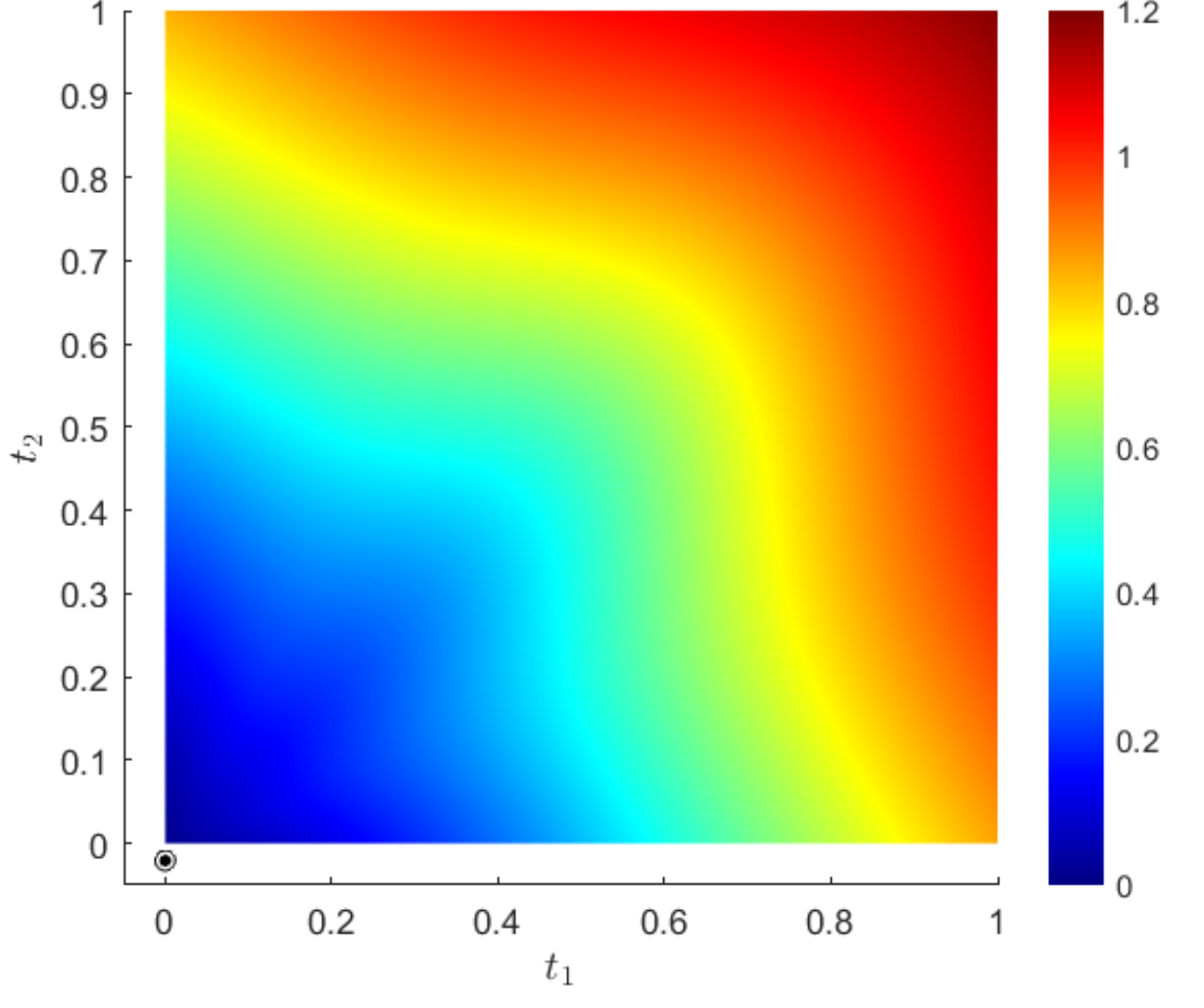}}
{\includegraphics[trim = 4cm 9.25cm 4.5cm 10cm, clip = true, height=5.5cm]{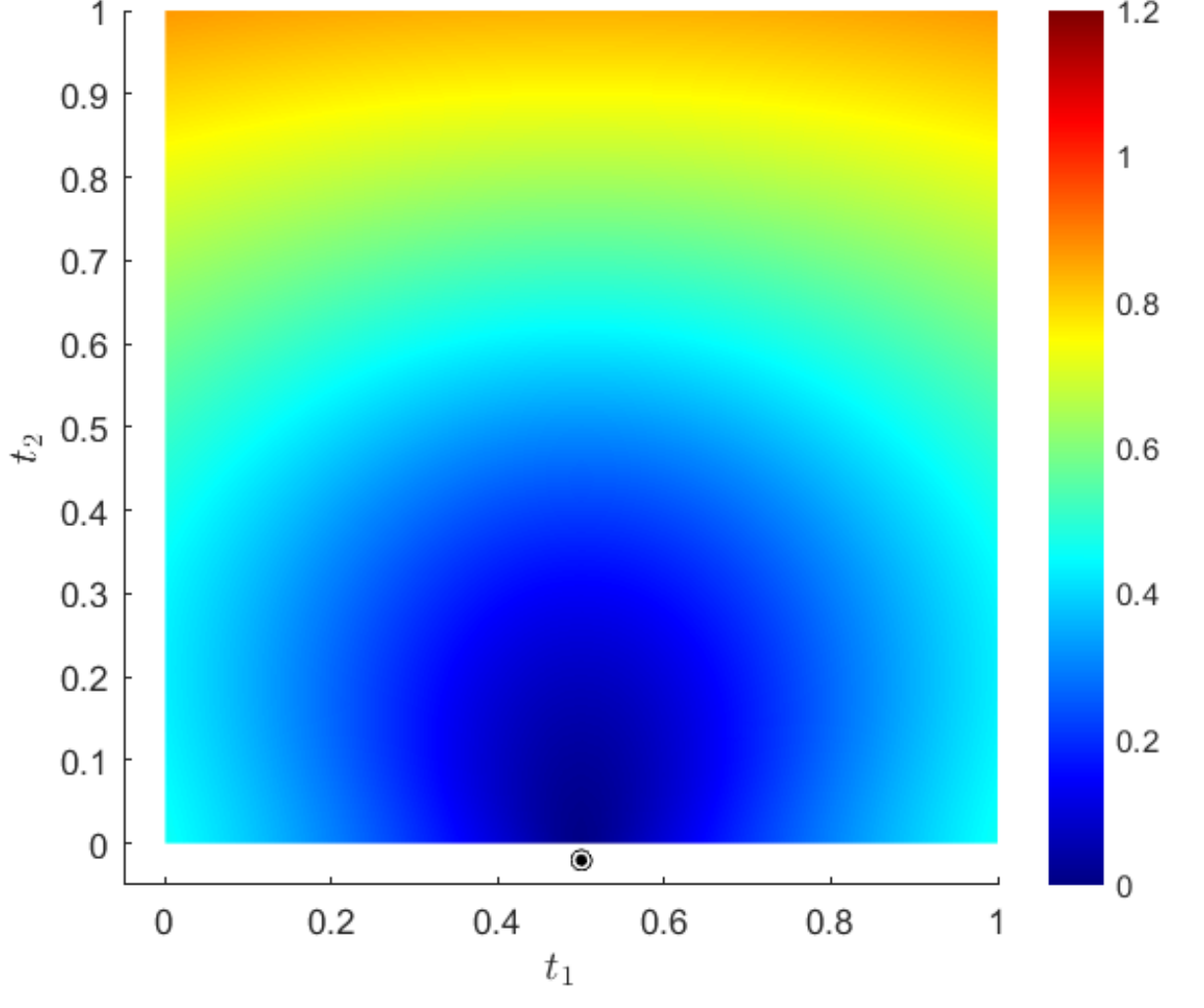}}
\caption{ An example of the absolute value of the regularized function of the nearly singular kernel with the source point and integration domain belonging to two different adjacent patches of a sphere. The top and bottom figures correspond to regularization for the single and double layer potential, respectively (with $m=2$). Two different configurations of points $\s_e$, represented as a black dots, are considered on the left and on right.}
\label{fig:extrapolation}
\end{figure}

\subsection{B-spline quasi-interpolation scheme}

To approximate integrals in 3D collocation BEM, we use the 2D Quasi Interpolation (QI) scheme firstly introduced in 
\cite{BGMS16} which is the tensor-product extension of the univariate Hermite QI scheme introduced in \cite{MSbit09}. With such approach a bivariate function $f$ is approximated on a rectangle $R$ by using a QI spline tensor product operator $  f \rightarrow \sigma_f := \sum_{ \i \in {\cal J}_R } \lambda_\i \hat B_{\i,\p}$, where $ \{ \hat B_{\i,\p} , \i \in { \cal J}_R \} $ is the tensor-product B-spline basis generating the chosen finite dimensional space of splines defined on $R,$ with ${ \cal J}_R$ denoting a local set of  multi--indices used to identify the basis elements.  In the considered QI scheme each coefficient $\lambda_\i, \i \in {\cal J}_R,$ is computed  as a  linear combination of  the values assumed  in a suitable local subset of breakpoints  by $f,$ by both its first partial derivatives and also by its second mixed derivative, see \cite{BGMS16} for details.  However, in the present work we rely on a derivative-free variant of such QI scheme, since we are interested to employ it for numerical integration. Such variant is obtained approximating the required derivative values  by suitable finite difference formulas. Then, collecting with lexicographical order all the coefficients $\lambda_\i, \i \in {\cal J}_R$ in a vector $\bfl$, in compact matrix notation we have
\begin{align}
\label{QIlam}
\bfl = \left( {\hat C}^{ (p) } \right)^\top  \f\,,
\end{align}
where $\f$ is a vector collecting with the same ordering all the values of $f$ at the breakpoints and
\begin{align}
\label{Cp}
{\hat C}^{(p)}  :=   \left( ({\hat A_2}^\top \otimes {\hat A_1}^\top)  -  ({\hat A_2}^\top \otimes ({\hat B_1}^\top D_{t_1})) -  (({\hat B_2}^\top D_{t_2}) \otimes {\hat A_1}^\top) +  (({\hat B_2}^\top D_{t_2}) \otimes ({\hat B_1}^\top D_{t_1}))  \right)^\top\,.  
\end{align}  
In the above formula $\hat {A_i},\, \hat {B_i}$ denote the two banded matrices, explicitly reported for example in \cite{CFSS18}, which are associated to the univariate original Hermite scheme formulated in the $i$-th direction, $i=1,2;$ $D_{t_i}$ are the matrices associated with the two analogous directional finite difference formulas used by the considered QI variant of the original scheme for approximating the first partial derivatives. Clearly, by incorporating a sufficiently accurate finite differences formula to approximate the required derivatives, the resulting derivative-free QI scheme shares  the approximation power with the original scheme.
On this concern, referring to the considered case with uniform breakpoints on $R$, the following error bound can be obtained for sufficiently smooth functions \cite{BGMS16}, 

\begin{align} \label{errorboundQI}
\Vert f \,-\, \sigma_f \Vert_{\infty,R} \le \,
\Vert \partial ^{p_1+1,0} f \Vert_{\infty,R} \left( \frac{H_1 }{\nu_1-1} \right)^{p_1+1} + \Vert \partial ^{0,p_2+1}f\Vert_{\infty,R} \left( \frac{H_2 }{\nu_2-1} \right)^{p_2+1}  \\ 
+\Vert  \partial^{p_1+1,p_2+1} f \Vert_{\infty,R} \left( \frac{H_1 }{\nu_1-1} \right)^{p_1+1} \left( \frac{H_2 }{\nu_2-1} \right )^{p_2+1} \,,\nonumber \end{align}
where $ H_i$ and $\nu_i$ respectively denote the size of $R$ and the number of uniform breakpoints in the $i$-th coordinate direction, $i=1,2,$  and $\partial^{i_1,i_2} f$ stands for ${ \partial^{i_1+i_2} f}/{(\partial t_1^{i_1} \partial t_2^{i_2})}.$ 
Thus the convergence order of the scheme with respect to $\max \{H_1 , H_2\} \rightarrow 0 $ but also to $\max \{ (\nu_1 - 1)^{-1}, (\nu_2 -1)^{-1}\}$ is equal to $\min \{p_1,p_2\}+1$.

\subsection{B-spline quadrature rules for regular integrals} 
\label{sec:regularRule}
 
Since in the assayed BIEs the regular integrands consist of the product between a bivariate regular function $f$ with a tensor product B-spline $\hat B_{\j,\d}$, we need to construct a quadrature rule to evaluate numerically the following integral: 
\begin{equation} \label{regint}
\int_{R_\j} f(\t) \hat B_{\j,\d}(\t)\,d\t.
\end{equation}

The first step consists in approximating $f$ with the QI scheme previously described. In this way we can write:
\begin{align*}
f \approx \sigma_f = \sum_{\i \in {\cal J}_{R_\i}}\lambda_\i \hat B_{\i,\p},
\end{align*}
where $\hat B_{\i,\p}$ denotes a tensor product B-spline of bi-degree $\p = (p_1,p_2)$, which is a trial function of a local spline space defined in  $R_\j$. 
Hence
\begin{align}\label{firstStep}
\int_{R_\j} f(\t) \hat B_{\j,\d}(\t)\,d\t \approx \int_{R_\j}\sigma_f(\t) \hat B_{\j,\d}(\t)\,d\t .
\end{align}
The product $\sigma_f \hat B_{\j,\d}$ is a new spline that can be expressed in terms of another B-spline basis defined in a suitable local ``product'' spline space of bi-degree $\p+\d$. Generalizing the idea developed in \cite{Morken91,falini2019adaptive} for the one-dimensional case, the coefficients of $\sigma_f\hat B_{\j,\d}$ in the B-spline basis of the product space can be written as $G_{\j}^\top \bfl,$ with $\bfl$ defined by QI as specified in (\ref{QIlam}). 
In the current bivariate setting the matrix $G_\j$ is defined as $G_\j := G_{j_1}^{(p_1,d_1)}\otimes G_{j_2}^{(p_2,d_2)},$ where $G_{j_1}^{(p_1,d_1)}$ and $G_{j_2}^{(p_2,d_2)}$ are suitable direction-wise coefficient matrices and $\j = (j_1,j_2)$. The dimension of  $G_\j$ is  $ \vert {\cal J}_{R_\j} \vert \times N_{\pi_\j}$, where $N_{\pi_\j}$ denotes the dimension of the introduced local product spline space. Thus, the integral in \eqref{firstStep} can finally be evaluated as,
\begin{align*}
\int_{R_{\j} } f(\t) \hat B_{\j,\d} (\t) d\t \approx  \vv^\top G_\j^\top \bfl =  \vv^\top G_\j^\top (\hat C^{(p)})^\top \f,
\end{align*}
where the matrix $\hat C^{(p)}$ is defined in (\ref{Cp}) and $\vv$ denotes a vector of suitable length collecting in lexicographical order the following integrals of each  B-spline $\hat B_{{\bf k},\p+\d}$ of the basis of the local product spline space,
\begin{align*}
\displaystyle{ \int_{R_\bp} \hat B_{{\bf k},\p+\d}(\t)\, d\t  \,=\, \frac{| \mbox{supp}(\hat B_{{\bf k},\p+\d}) |}{(p_1+d_1+1) (p_2+d_2+1)} }\,, \qquad \k \in {\cal J}_{\pi_\j} \,,
\end{align*}
where ${\cal J}_{\pi_\j}$ is a set of $N_{\pi_j}$ distinct multi-indices, as usually done to identify a tensor-product B-spline basis, and ${R_\bp}$ is the support of the spline.
 
As already mentioned, the derived integration rule is the tensor product extension of a variant specific for integrals including a B-spline weight of the 1D quadrature formula introduced in \cite{MS12}. Therefore, the error analysis developed in \cite{MS12} can be easily  extended to such variant and lifted to the tensor product setting. To be more concise, we assume that a uniform distribution of quadrature nodes is used in each $R_\j$ (i.e., uniform breakpoints in the associated local spline space are used by the QI scheme), since this is our common setting in numerical experiments. Assuming $f$ is a sufficiently smooth function, the quadrature error of the rule to approximate the integral in \eqref{regint} is upper bounded by a positive constant times the following quantity,
\begin{align} \label{errorbound}
\frac{ H_{1,\j} H_{2,\j} }{(d_1+1)(d_2+1)}\, \left\{\, \Vert \partial ^{p_1+r_1+1,0} f\Vert_{\infty,R_\j} \left( \frac{H_{1,\j} }{\nu_1-1} \right)^{p_1+r_1+1} + \Vert \partial ^{0,p_2+r_2+1}f\Vert_{\infty,R_\j} \left( \frac{H_{2,\j} }{\nu_2-1} \right)^{p_2+r_2+1} \right. \\ 
+ \left. \Vert  \partial^{p_1+r_1+1,p_2+r_2+1} f \Vert_{\infty,R_\j} \left( \frac{H_{1,\j} }{\nu_1-1} \right)^{p_1+r_1+1} \left( \frac{H_{2,\j} }{\nu_2-1} \right )^{p_2+r_2+1} \right\}\,,\nonumber \end{align}
where $ H_{i,\j}$ and $\nu_i$ respectively denote the size of $R_{\j}$ and the number of uniform quadrature nodes in the $i$-th coordinate direction, $i=1,2$ and $r_i=1 (0)$ for $p_i$ even (odd). 

For fixed numbers $\nu_i$ of nodes, the convergence order of the quadrature scheme with respect to $\max \{H_{1,\j} , H_{2,\j}\} \rightarrow 0 $ is 
\begin{align*}
C_1 \min \{H_{1,\j}^{p_1+r_1+3}, H_{1,\j}^{p_2+r_2+3} \},
\end{align*}
if the quotient $H_{1,\j}/ H_{2,\j}$ is bounded from below and above for all nested integration subdomains $R_\j$. The constant $C_1$ contains norms of all the involved derivatives of $f$ in \eqref{errorbound} for the largest considered $R_\j$. For fixed $R_\j$, the convergence order with respect to $\min \{ \nu_1, \nu_2 \}$ is equal to 
\begin{align*}
\tilde C_1 \min \{ (\nu_{1}-1)^{-p_1-r_1-1}, (\nu_{2}-1)^{-p_2-r_2-1} \},
\end{align*}
where the positive constant $\tilde C_1$ contains also constants $H_{i,\j}$ but not varying $\nu_i$.

\subsection{B-spline quadrature rules for singular integrals}
\label{sec:regularSing}
Let us focus now on the rules developed for the numerical approximation of the second addend on the right of \eqref{eq:subtraction}, which involves the simplified singular kernel $U_\s^m(\s-\t).$
The first step is analogous to that adopted for regular integrals, since $g$ is approximated with $\sigma_g$ by using the same QI approach and then multiplied by $\hat B_{\j,\d}.$ Thus we can write
\begin{align} \label{simplsing}
\int_{R_\j}{U}_{\s}^m(\s-\t) \hat B_{\j,\d} (\t) g(\t)\,d\t \approx \int_{R_\j}{U}_{\s}^m(\s-\t) \hat B_{\j,\d} (\t) \sigma_g(\t)\,d\t  = \vv^\top G_\j^\top (\hat C^{(p)})^\top \g\,,
\end{align}
where $\g$ is the vector whose entries are the values of $g$ used by the QI operator but now $\vv$ collects in a suitable order all the following  integrals
\begin{align} \label{moments}
\int_{R_\j}{U}_{\s}^m(\s-\t) \hat B_{\k,\p+\d}(\t)\ d\t\,, \qquad \k \in {\cal J}_{\pi_\j}.
\end{align}
Thus, in order to complete the definition of the quadrature rule, we are interested in deriving the exact expression of the above integrals \eqref{moments}. For this aim we introduce the following notation, 
$$ I^{q_1,q_2}_{r_1,r_2}(k_1,k_2) := \int_{\tau^{(1)}_{k_1}}^{\tau^{(1)}_{k_1+r_1+1}} \int_{\tau^{(2)}_{k_2} }^{\tau^{(2)}_{k_2+r_2+1}}{U}_{\s}^m(\s-\t) (t_1-s_1)^{q_1}(t_2-s_2)^{q_2} \hat B_{k_1,r_1}(t_1)\hat B_{k_2,r_2}(t_2) dt_1 dt_2\,,$$
where $\k = (k_1,k_2)$ and $\tau^{(\ell)}_{k_\ell+j}, j=0,\ldots,r_\ell, \ell=1,2$ define the two univariate knot vectors active in the definition of $\hat B_{\k,\p+\d},$ where for brevity we have set $r_\ell := p_\ell+d_\ell $, which is the univariate degree  in the $\ell$-th direction of the product space. Then, using the tensor product extension of the Cox-de Boor recurrence relation for B-splines \cite{deBoor01},
we can express $I^{q_1,q_2}_{r_1,r_2}(k_1,k_2)$ as a linear combination of $I^{q_1+m_1,q_2+m_2}_{r_1-1,r_2-1}(k_1+w_1,k_2+w_2)$, where $m_1,m_2,w_1,w_2=0,1$.
Starting with $q_1=q_2=0$ and iterating this recurrence up to B-splines of degree 0 in both coordinate directions, the procedure can be completed, provided that the following {\it basic moments} have been preliminarily computed for $q_\ell$  and $i_\ell$, respectively ranging in the set of indices $ 0,\ldots, p_\ell+r_\ell$ and $k_\ell,\ldots,k_\ell+r_\ell, \ell=1,2,$
\begin{align*}
I^{q_1,q_2}_{0,0}(i_1,i_2) &=  \int_{\tau^{(1)}_{i_1}}^{\tau^{(1)}_{i_1+1}} \int_{\tau^{(2)}_{i_2} }^{\tau^{(2)}_{i_2+1}}{U}_{\s}^m(\s-\t) (t_1-s_1)^{q_1}(t_2-s_2)^{q_2}  dt_1 dt_2 \\
\ &=  \int_{s_1-\tau^{(1)}_{i_1+1}}^{s_1-\tau^{(1)}_{i_1}} \int_{s_2-\tau^{(2)}_{i_2+1}}^{s_2-\tau^{(2)}_{i_2} }{U}_{\s}^m(\z) (-z_1)^{q_1}(-z_2)^{q_2}  dz_1 dz_2.
\end{align*}
For this aim we rely on the analytical expressions derived in \cite{tadej21} computed with the help of Wolfram Mathematica software.

We are left to derive the asymptotic accuracy of these rules with respect to the $\max\{H_{1,\j}\,,H_{2,\j}\} \rightarrow 0,$ where, as in the previous subsection, $H_{i,\j}$ denotes the size of the integration domain $R_\j$ in the $i$-th coordinate direction. Now from (\ref{simplsing}) we can easily derive that the absolute value of the related quadrature error can be upper bounded by 
\begin{align}\label{err_sing}
\Vert g - \sigma_g \Vert_{\infty,R_\j} \ \int_{R_\j} \vert U_{\s}^m(\s-\t) \vert \hat B_{\j,\d} (\t) \,d\t.
\end{align}
Using the triangle inequality $\vert U_{\s}^m(\s-\t) \vert \le \vert U(\s,\t) \vert + \vert  U_{\s}^m(\s-\t) - U(\s,\t) \vert$, we can bound the integral in \eqref{err_sing} by a sum of two new integrals. Then, by applying  Proposition~\ref{Prop2} reported in the Appendix, we can estimate the first integral as
\begin{align}\label{est1}
 \ \int_{R_\j} \vert U(\s,\t) \vert \hat B_{\j,\d} (\t) \,d\t \leq  
 \ \int_{R_\j} \vert U(\s,\t) \vert  \,d\t \leq
 C_2 \max\{H_1,H_2\}.
 \end{align}
The second integral involves a tail of the series expansion of $U(\s,\bullet)$ about the source point $\s$. Referring to results in \cite{tadej21}, on this concern we can state that there exists a positive constant $C_3$ involving the infinity norm of derivatives of $\tilde \F$ of degree $\leq m$ on $R_\j$, such that
\begin{align*}
\Vert U_{\s}^m(\s-\t) - U(\s,\t) \Vert_{\infty, R_\j} \le C_3 \frac{ \left(\max\{H_{1,\j}\,,\, H_{2,\j}\} \right)^{m-1}}{m!},
 \end{align*}
 hence 
\begin{align}\label{est2} 
\int_{R_\j} \vert  U_{\s}^m(\s-\t) - U(\s,\t) \vert \hat B_{\j,\d} (\t) \,d\t \leq C_3 \frac{ \left(\max\{H_{1,\j}\,,\, H_{2,\j}\} \right)^{m}}{m! (d_1 +1)( d_2 + 1)}.
\end{align}

By combining these estimates with the bound in (\ref{errorboundQI}) for the QI error, we finally derive the following error bound for the singular integral in \eqref{simplsing} with respect to size of $R_\j$,
\begin{align*}
C \max \{H_{1,\j}^{p_1+2}, H_{2,\j}^{p_2+2}\},
\end{align*}
where $C = C_1 \max\{C_2,C_3\}$ and $C_1$ contains norms of all the involved derivatives of $f= g - \sigma_g$ in \eqref{errorbound} for the largest considered $R_\j$. Instead, by fixing domain $R_\j$, the error bound with respect to $\min \{ \nu_1-1, \nu_2 -1\}$ is equal to
\begin{align*}
\tilde C \max\{ (\nu_1-1)^{-p_1-1}, (\nu_2-1)^{-p_2-1} \},
\end{align*}
where $\tilde C$ contains a corresponding constant $\tilde C_1$ to the one at the end of Section~\ref{sec:regularRule}, involving the function $f= g - \sigma_g$, multiplied by right-hand side estimates of \eqref{est1} and \eqref{est2}.


\subsection{The threshold for near singularity detection} \label{thresh}
The efficient and accurate evaluation of different types of integrals \eqref{genint}  appearing in the considered BIEs is crucial, see also \cite{sorgentone22} on such concern. For this aim in this subsection we are interested in determining a reasonable threshold $\delta$ to be used in (\ref{eqn:thresh_dist}) to switch between nearly singular and regular integrals. 
The goal is to apply  the error bound \eqref{errorbound}  introduced in Section~\ref{sec:regularRule} to $f=gU,$ in order to derive on each patch a mesh dependent definition of the threshold $\delta$ to be used in (\ref{eqn:thresh_dist}).

By analyzing the nearly singular nature of the considered integral in \eqref{genint}, it is important to observe that the function $f:=gU (\s, \bullet)$ for a fixed $\s=\tilde \F^{-1}(\tilde \x)$ includes a factor proportional to $1/\tilde r$, where $\tilde r = \| \tilde \x -\tilde \F(\bullet) \|_2$. Hence, $\Vert  \partial ^{i_1,i_2}f \Vert_{\infty,R_\j}$ includes a term $\tilde r_{\rm{min},\j}(\x)^{-i_1-i_2-1},$ where $\tilde r_{\rm{min},\j}(\x)$ is the minimal distance between $\tilde \x$ and the physical integration domain image of $R_\j$ via the geometry map  $\tilde \F,$  see its definition in formula \eqref{eqn:thresh_dist}. Thus, considering the expression of the error bound in \eqref{errorbound} for a fixed number of uniformly spaced quadrature nodes, we can infer a remarkable reduction of the convergence order with respect to $\max \{H_{1,\j},H_{2,\j} \} \rightarrow 0$ if $\tilde r_{\rm{min},\j}(\x)$ would reduce proportionally to it, due to increasing terms ``$\Vert \partial ^{i_1,i_2} f \Vert_{\infty,R_\j}$''. Instead, by setting
\begin{align*}
\tilde r_{{\rm min},\j}(\x) \ge  \eta_\j:= \max \left\{ H_{1,\j}^{\frac{1}{2(p_1+r_1+2)}}\,,\, H_{2,\j}^{\frac{1}{2(p_2+r_2+2)}} \right\}
\end{align*}
this guarantees that for all the scaled integrals treated with the regular rule in the governing BIEs, the quadrature convergence order is at least equal to $5/2+ \min \{ p_1+r_1, p_2+r_2 \}.$
In order to have just one threshold $\delta$ for all the integrals involving the same  geometry map, for the $k$-th patch we set
\begin{align*}
\delta = \delta_k := \max_{\j \in {\cal I}_k}  \eta_\j.
\end{align*}

\section{Numerical examples} \label{sec:NE}


In this section we test our model to numerically solve one interior and two exterior Helmholtz problems, all with Neumann boundary conditions, which are more commonly considered for applications. 

The two well-known benchmarks  for exterior problems are the pulsating sphere and the rigid scattering on a sphere; see for instance \cite{Simpson14,Venas20,Anitescu21}. 
The  sphere has radius $R=1$ and it is parameterized by $6$ quartic ($\d_g = (4,4)$) NURBS patches, see Figure \ref{fig:patch}(a) and refer to 
\cite{Cobb1994TilingTS} for details.  As in \cite{Dolz18}, the considered quartic NURBS parameterization based on cube topology is the preferred choice to the 8-patch tiling in order to avoid singularities in the geometry description.
The last example is a standard acoustic problem interior to a torus; see for example \cite{Simpson14,Venas20}. The toroidal surface $\Gamma$ has inner and outer radii respectively equal to $1$ and  $3$ and is exactly represented in a $16$-patch quadratic ($\d_g = (2,2)$) NURBS form, see Figure~\ref{fig:patch}(b) and refer to \cite{roma} for details.
\begin{figure}[t!]
\centering
\subfigure[The $6$-patch quartic NURBS representation of a sphere.]{
\includegraphics[trim = 0.0cm 0cm .0cm .5cm, clip = true, height=4cm]{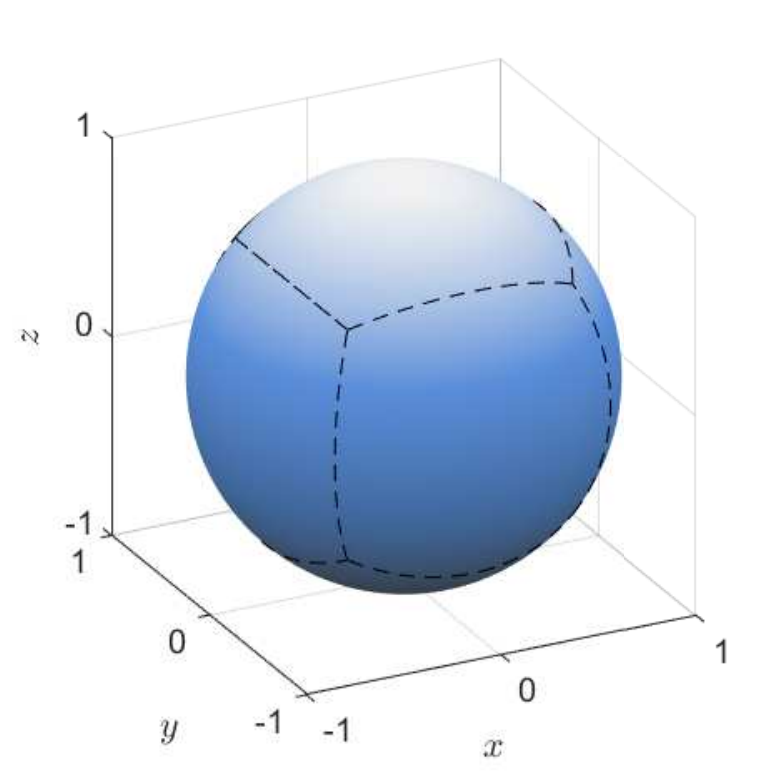}
}
\qquad
 \subfigure[The $16$-patch quadratic NURBS representation of a torus.]{
\includegraphics[trim = 0.0cm 0.5cm 0cm .5cm, clip = true, height=4cm]{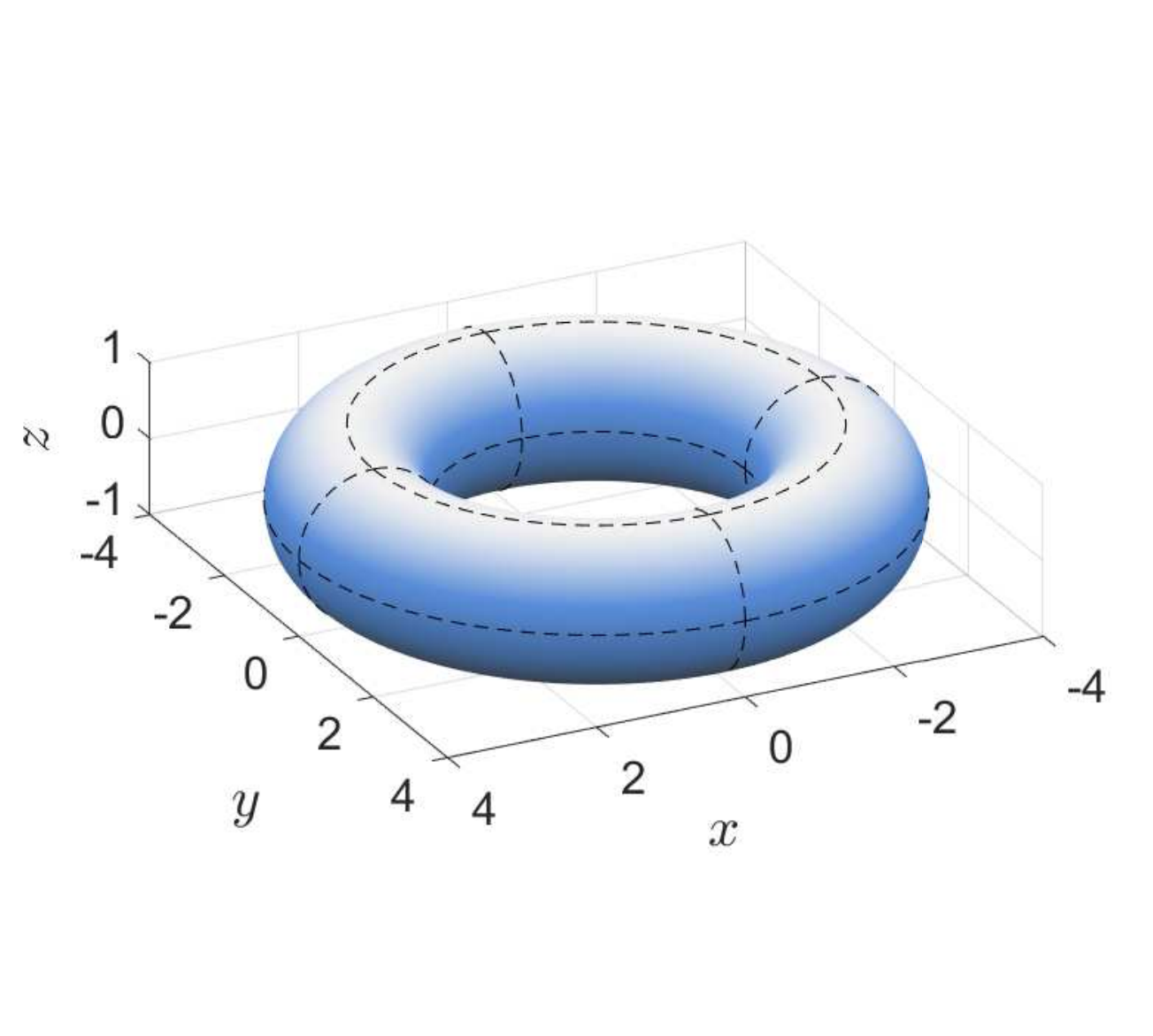}
}
\caption{The two NURBS multipatch representations of a sphere and a torus.}
\label{fig:patch}
\end{figure}


For the discretization space $V$, on each patch we employ standard tensor product B-splines of different bi-degrees $(d ,d)$ and of highest continuity; in the last two experiments we enforce $C^0$ continuity for the basis functions across patch interfaces. 
The quality of the numerical results is evaluated by considering the behavior of the on-surface relative error measured in the $L^2(\Gamma)$ norm,
 $$e_{L^2} = \frac{\Vert u_h - u^{ex} \Vert_{L^2(\Gamma)}}{\Vert u^{ex} \Vert_{L^2(\Gamma)}}.$$



\subsection{Pulsating sphere}

 The following problem serves as a benchmark test for the developed quadrature rules; it is also referred to as a patch test for IgA-BEM. The model can be used for example to numerically compute the sound pressure at a distance $r$ from the center of the sphere for a constant speed of wave in the radial direction. For the Neumann boundary conditions set to $u_{\rm N} \equiv e^{i \kappa} (i \kappa -1) /(4  \pi)$ for $ \x \in \Gamma$, the analytic solution is given
as $u = e^{i \kappa r} /(4 \pi r)$, hence the missing Cauchy data is $\phi \equiv  e^{i \kappa} /(4 \pi)$ for $ \x \in \Gamma$.

Since the geometry representation of the sphere is exact and $\phi$ is constant over the entire boundary, the error can be attributed solely to numerical integration. On each patch the knot vectors describing the local part of the discretization space $\hat {\cal S}_{h^{(\ell)}}^{(\ell)}$ are
\begin{align*}
T_{1}^{(\ell)} = T_{2}^{(\ell)} = \begin{bmatrix}0 & 0.25 & 0.5 & 0.75 & 1 \end{bmatrix}.
\end{align*}
The collocation points are midpoints of the knots, mapped to the physical domain via geometry map. Fixing $\kappa=1, \d = (0,0)$ and setting $c = 0.25$ for completing the definition of the threshold used for near singularity detection, we measure the $L^2$ error of the approximate solution with respect to different number of quadrature nodes on the support of B-splines. The number $n_q$ of quadrature nodes increases  in each coordinate direction with $ n_q = 12 \alpha+1 $, for $\alpha=1,2,\dots,6$. The same quadrature nodes are used for singular, regularized and regular integrals.

In Figure~\ref{fig:pulsatingSphere}(a) we observe the error distribution (modulus of the difference between the exact and approximate solution) for $\alpha=1$, $m=2$ (number of terms in the singularity extraction) , $p=2$ (spline degree of the QI operator). The error is well evenly spread on the boundary with a slight increase around patch corners, where the kernel expansion is more critical. The error is also slightly higher in the interior of the patches than in the neighborhood of the edges since the area of interior mapped cells is slightly larger than the boundary ones.
In Figure~\ref{fig:pulsatingSphere}(b) we observe that error convergence order is impacted by the choice of $m$ and $p$. For $p=m=2$ the obtained order of convergence is 3. If we increase both $m$ and $p$ to 3, the order of convergence increases to 4.

 \begin{figure}[thb]
\centering

\subfigure[Error distribution on the boundary for $\alpha=6, m=2,p=2$.]{
\includegraphics[trim = 0.25cm 0cm .0cm 0cm, clip = true	, height=5cm]{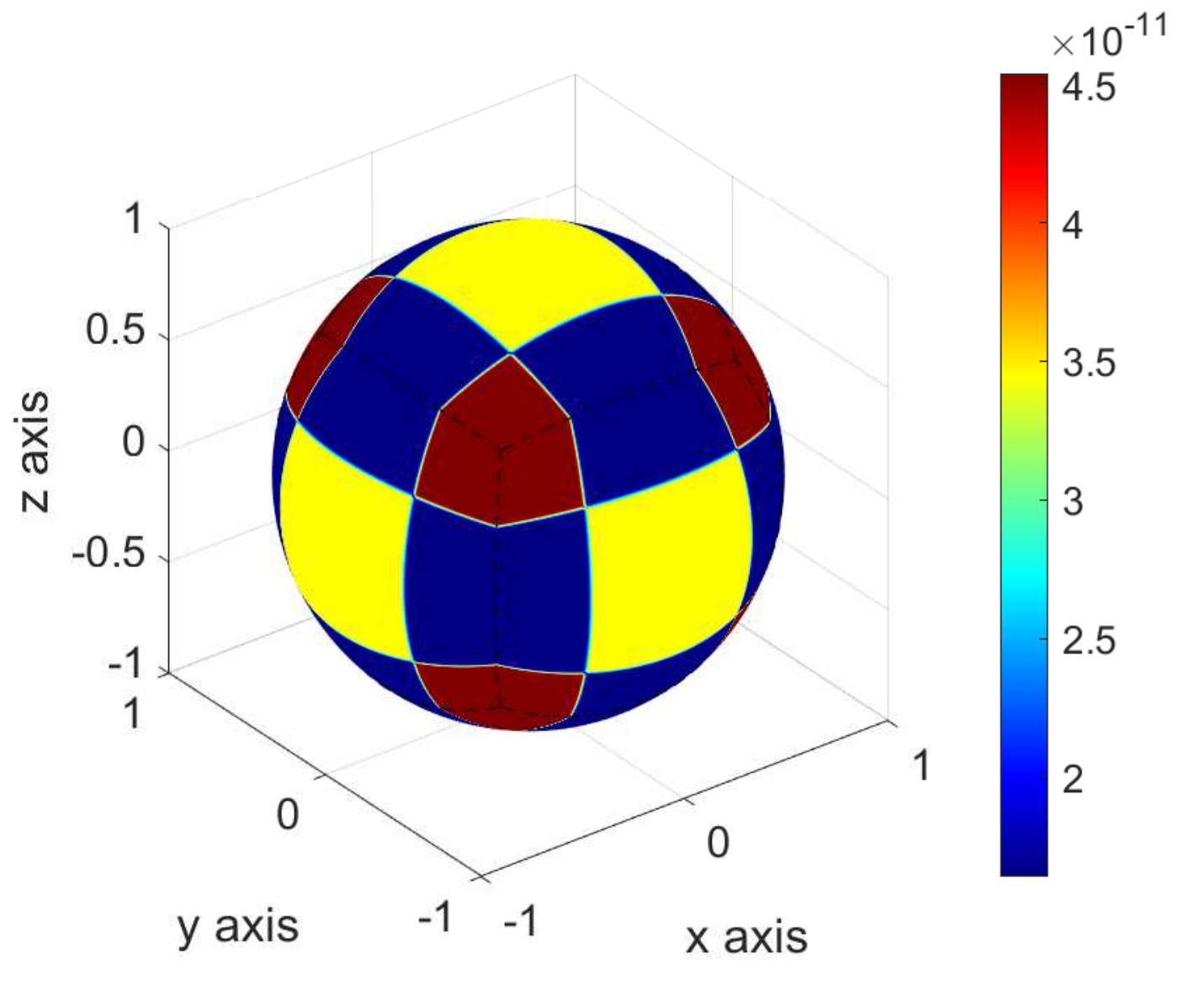}}
\qquad
 \subfigure[Error convergence plot with respect to multiplier $\alpha$ for quadrature number.]{
\includegraphics[trim = 0.0cm 0cm .0cm 0cm, clip = true, height=5cm]{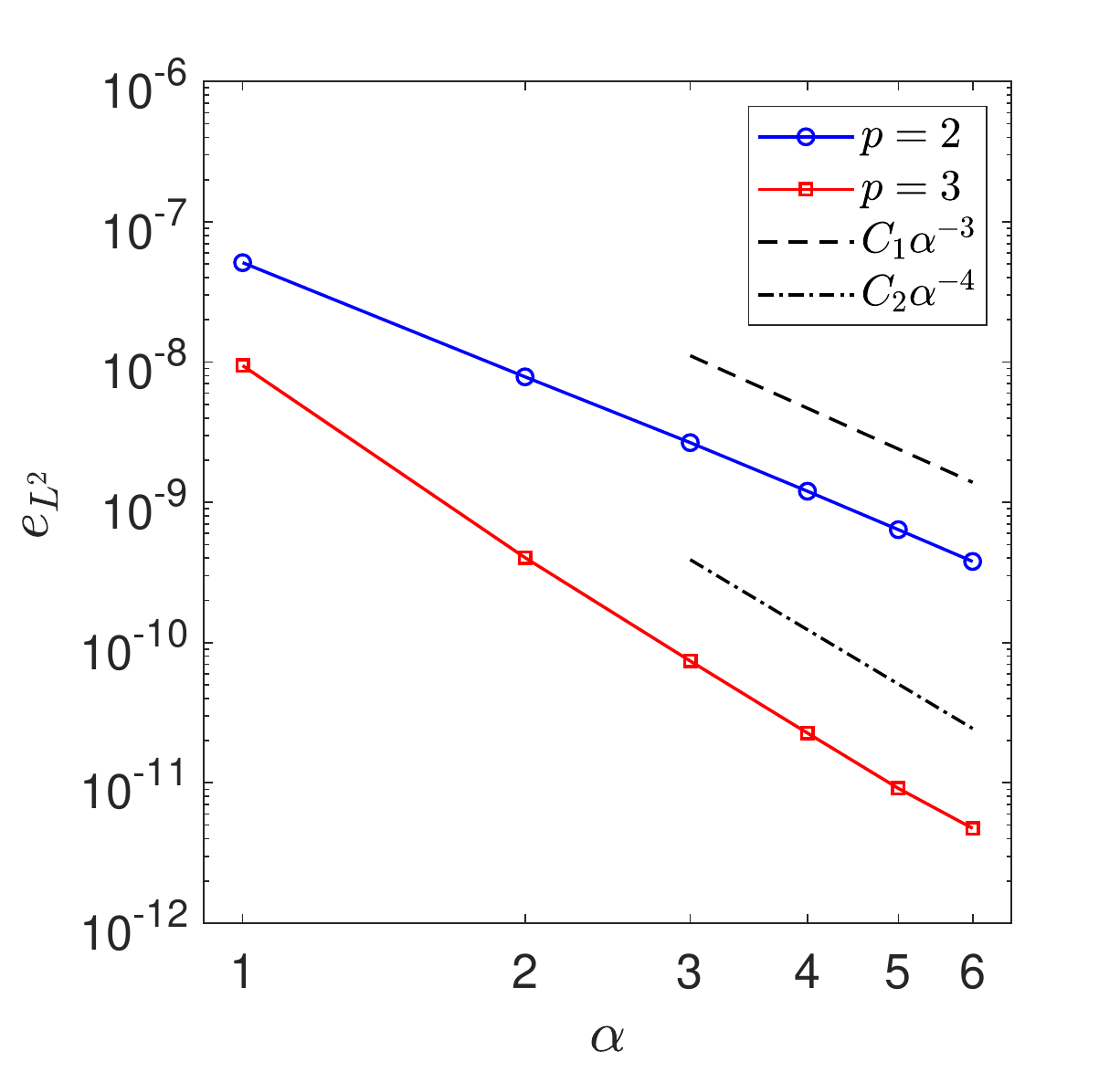}}
\quad\quad
\caption{Pulsating sphere for $\kappa=1$.}
\label{fig:pulsatingSphere}
\end{figure} 

The comparison can be done with Figure 6(a)  of \cite{Venas20}, where $R \kappa=1$ as we have used in our experiments. In the referred picture the results obtained with new approach (termed as {\it New}) are compared with those produced by the method introduced in \cite{Simpson14} (termed as {\it Old}). The relative $L^2$ error, $e_{L^2}$, is shown against the total number of quadrature nodes $n_{qp}$ (where the -negligible- number of quadrature nodes in the elements with active singularity is not taken into account).
In our case the total number of quadrature nodes, for $\alpha=6$, is $6\ \cdot 16\ \cdot 72^2\approx 5.0 \ \cdot 10^5$ achieving an error (with $p=3$) of  order $10^{-11}.$ This is a much better result than that produced by the {\it Old} approach but also better than that given by the {\it New} one.

\subsection{Rigid Scattering on a sphere}
In rigid scattering (exterior) problems reformulated in the frequency domain, it is assumed that an acoustic pressure $p_{inc}$ of amplitude $A$ is produced by a wave vector  $\kappa \vv,$ with $\vv$ denoting a unit vector prescribing the direction of the wave,
$$p_{inc}(\x) := A e^{i \kappa(\vv \cdot \x)}.$$
Such pressure is incident on a rigid body represented by a volume $\Omega^{(i)}$ with boundary $\Gamma,$ where the rigidity assumption implies that it reacts producing in $\Omega^{(e)}$ an additional scattered pressure field $p$ verifying the Helmholtz equation in $\Omega^{(e)}\,,$ the Sommerfield radiation condition at infinity and having a variation in the normal direction on $\Gamma$ opposite to that of $p_{inc}$, see for example \cite{Venas20}. The scattered pressure $p$ is the unknown to be determined in $\Omega^{(e)}$ and it can be written as the difference between a total pressure field $p_{tot}$ and $p_{inc},$ where $p_{tot}$ is such that
$$ \left\{ \begin{array}{ll} p_{tot}(\x) &= p_{inc}(\x) +  (V_\kappa \partial_n p_{tot}) - (K_\kappa p_{tot})(\x)\,,\,\, \x \in \Omega^{(e)} \cr
\partial_n p_{tot} (\x) &= 0\,,\,\,  \forall \x \in \Gamma\,. \cr
\end{array}
\right.$$
 Note that in this experiment we assume  $\Gamma$ as a sphere with radius $R=1$ (again parameterized as previously described) and so, for symmetry, without loss of generality we can assume $\vv = (1,0,0).$  The problem can be used as a benchmark because the analytic expression of $p$ is a priori known, see for example \cite{Anitescu21},
 $$ p(\x) = - A\ \sum_{n=0}^{\infty}\frac{i^n (2n+1) j_n'(\kappa R)}{h_n'(\kappa R)} \ P_n(\cos(\theta) )\ h_n(\kappa r)\,,$$
 where in our experiments we set $A=1$ and $n$ ranges from $1$ to $10$ because this is sufficient to study the convergence of our numerical scheme for the used frequencies. In this formula, assuming the sphere centered in   ${\bf O}$,  $ r  = \Vert \x - {\bf O} \Vert_2,$  $\theta$ is the angle between the vector $\x - {\bf O}$ and $\vv,$  $h_n$ and $j_n$ are the spherical Hankel and Bessel functions of order $n$ (the derivative with respect to its argument is denoted by $'$) and $P_n$ is the Legendre polynomial of order $n.$
 \begin{figure}[thb]
\centering
 \subfigure[The real part of the total pressure on the sphere]{
\includegraphics[trim = 0.6cm 0cm .75cm 0cm, clip = true, height=2.8cm]{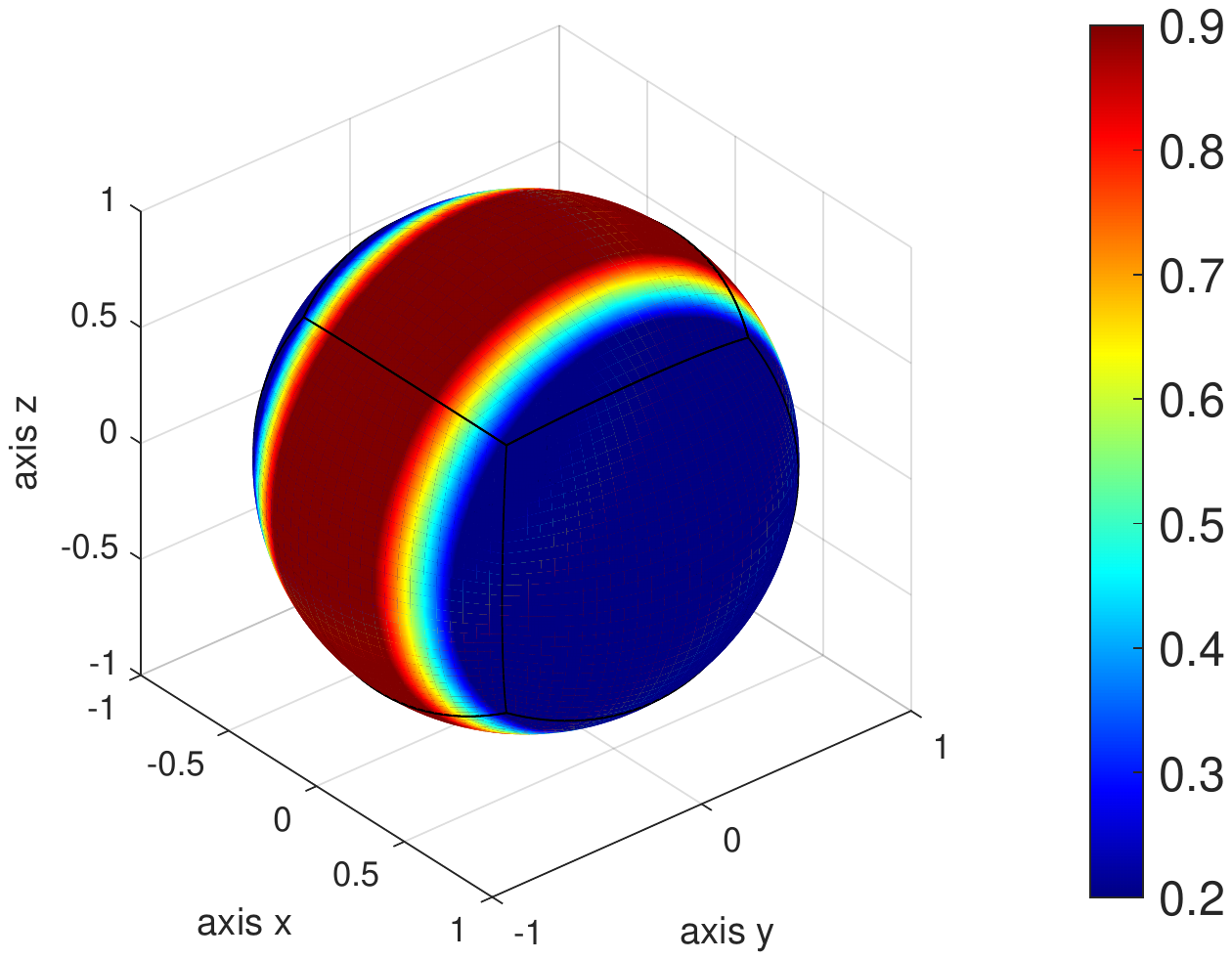}
}
\hspace{0.2cm}
\subfigure[The imaginary part of the total pressure on the sphere]{
\includegraphics[trim = 0.6cm 0cm .75cm 0cm, clip = true, height=2.8cm]{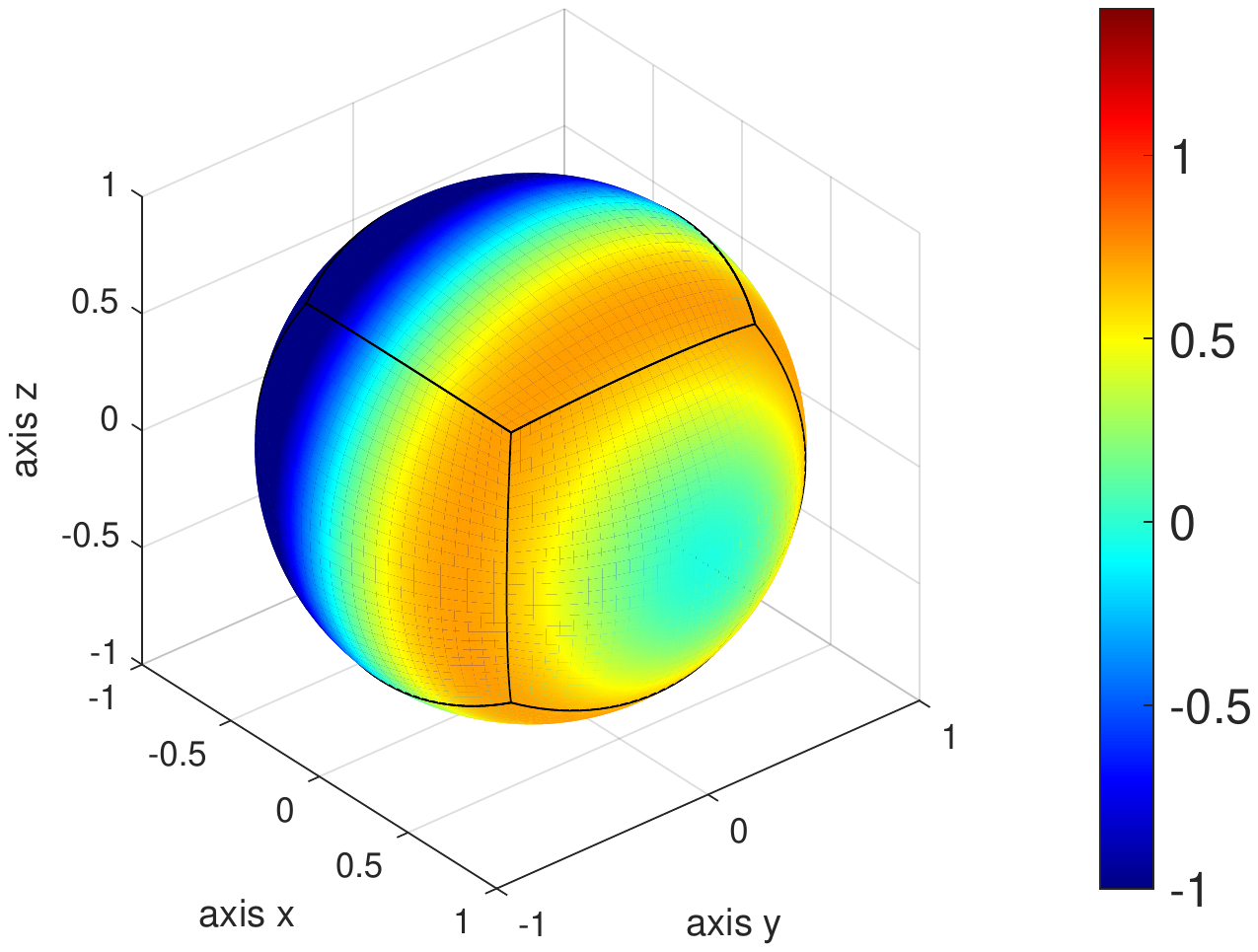}
}
\hspace{0.3cm}
 \subfigure[The real part of the total pressure on the equatorial plane]{
\includegraphics[trim = 0.6cm 0cm .75cm 0cm, clip = true, height=2.8cm]{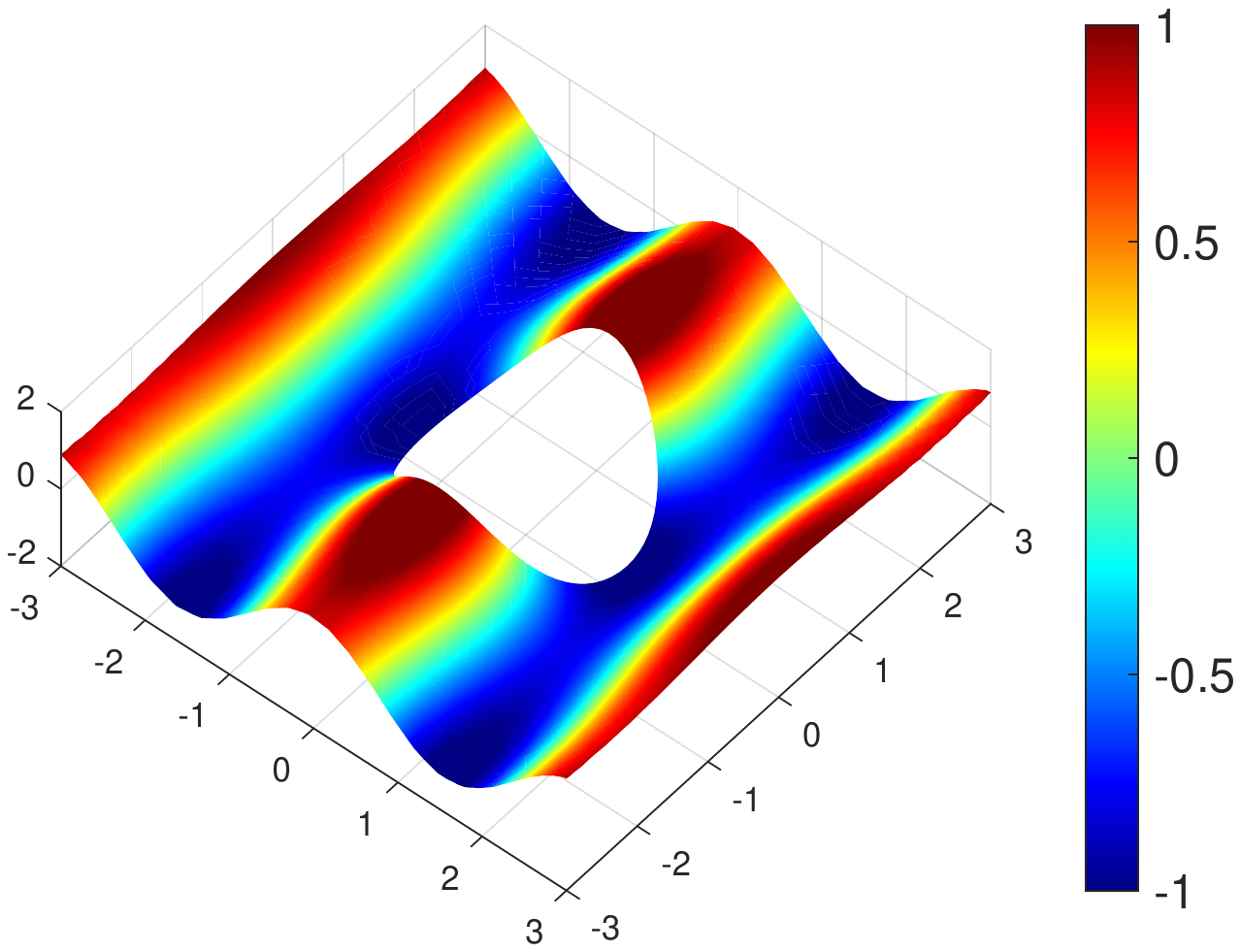}
}
\hspace{0.2cm}
 \subfigure[The imaginary part of the total pressure on the equatorial plane]{
\includegraphics[trim = 0.6cm 0cm .75cm 0cm, clip = true, height=2.8cm]{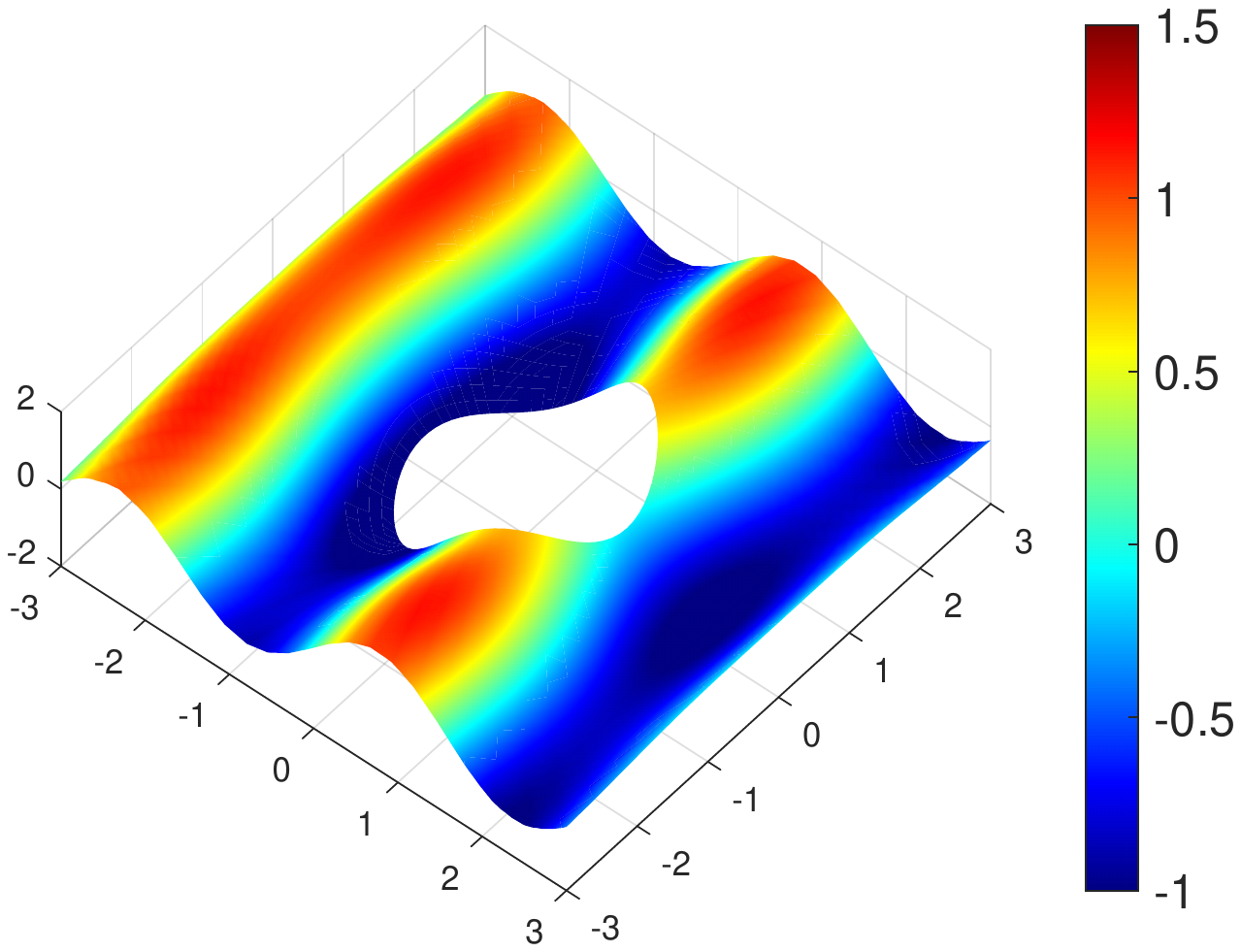}
}
\caption{Rigid scattering on a sphere with $\kappa=2$: on and off surface total pressure distributions.}
\label{fig:SS_sol}
\end{figure}
\begin{figure}[t!]
\centering
\vspace{2cm}
\subfigure[On surface error distribution ($d= 2\,, n=5$)]{
\includegraphics[trim = 0.6cm 0cm .75cm 0cm, clip = true, height=5cm]{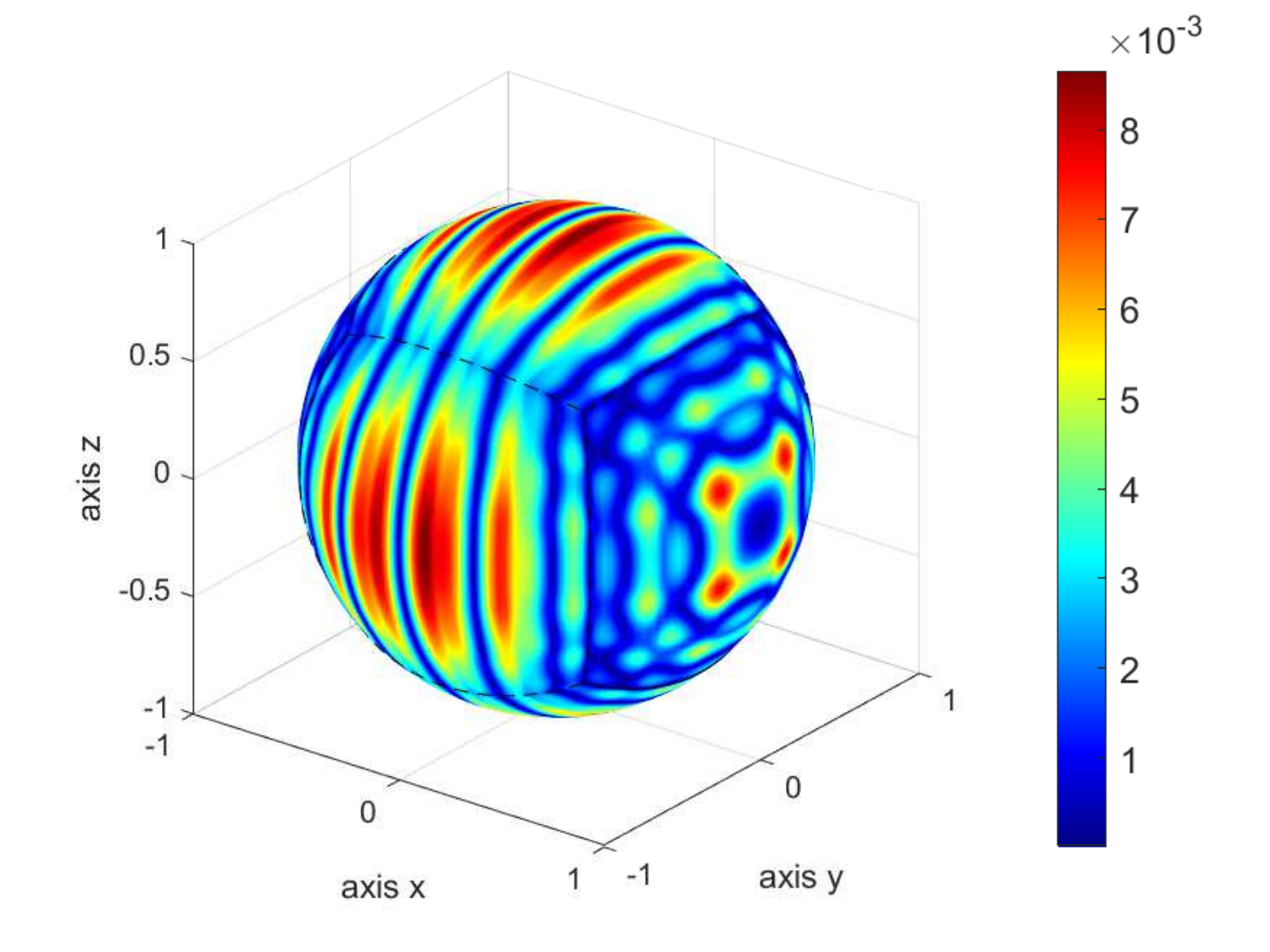}
}
\hspace{2cm}
\subfigure[On surface error distribution ($d= 4\,, n=3$)]{
\includegraphics[trim = 0.6cm 0cm .75cm 0cm, clip = true, height=5cm]{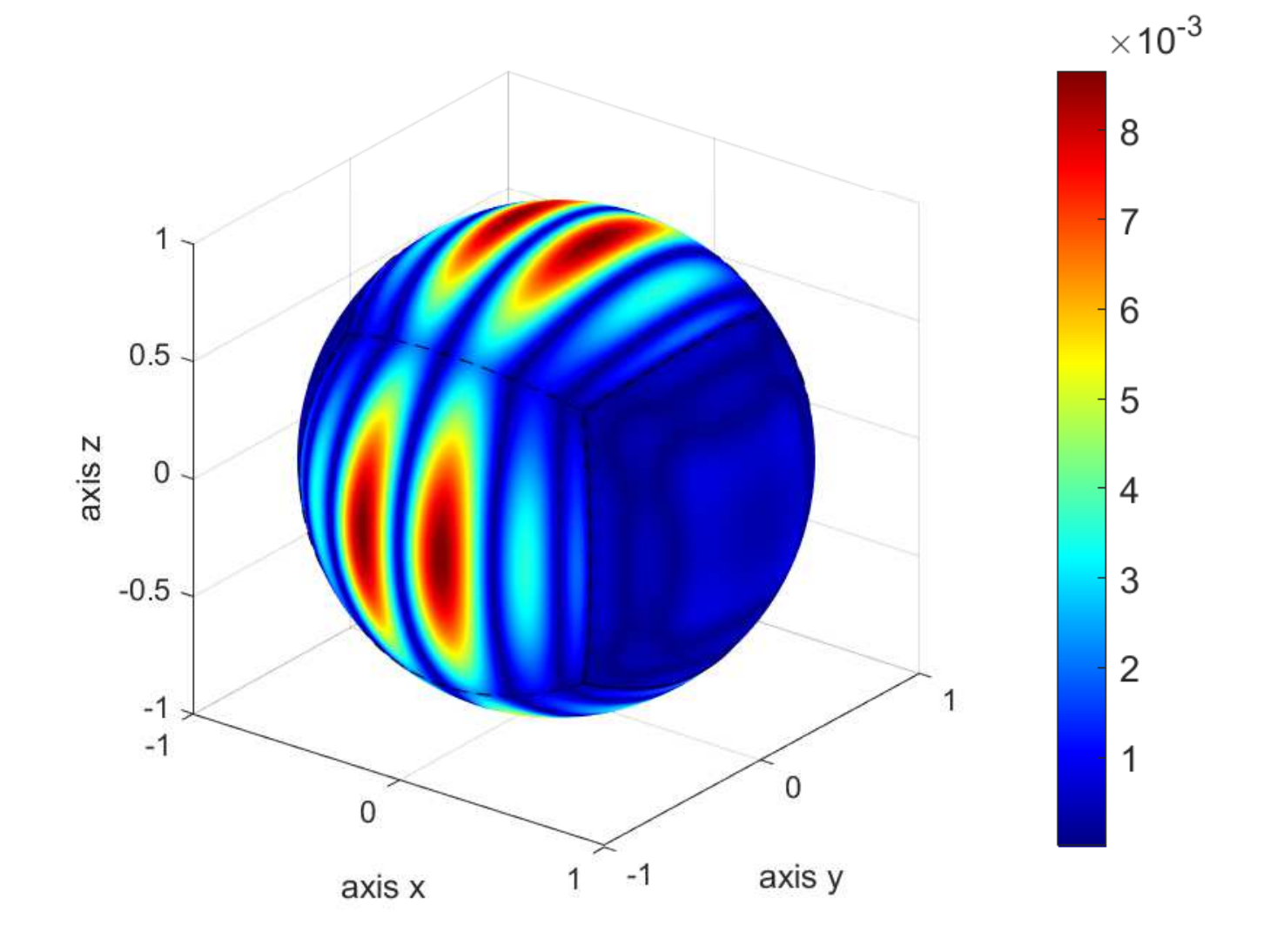}
}

 \caption{Rigid scattering on a sphere with $\kappa=2$: on surface distribution of the modulus of the absolute error  ($N_{\rm DOF}= 218$).  }
\label{fig:SS_error}
\end{figure}

%
The values considered for $\kappa$ in the experiments are $ \kappa= 1,2,3$ which in the air corresponds to a frequency of about $50, 100, 150$ Hz. In the implementation we aim to approximate the unknown Dirichlet  Cauchy datum $\phi = p_{tot}|\Gamma$ (alternatively also denoted below just with $u$ which denotes also its extension to $\Omega^{(e)}$). Thus, considering that the corresponding given Neumann datum is homogeneous,  the boundary integral equation to be considered simplifies to the following one,
\begin{equation} \label{rigidscat}
 (K_\kappa + \frac{1}{2}I) \phi (\x) \,=\, p_{inc}(\x)\,, \quad \x \in \Gamma\,.
 \end{equation}
The discretization of this BIE has been done by using a globally $C^0$ multi-patch IgA space with the same size of uniform elements on each patch. For this experiment we consider the following simplification
\begin{align*} 
\frac{\r \cdot \n_\y}{r^2}  = - \frac{1}{2R},
\end{align*}
where $\n$ points toward the center of the sphere since the considered problems are exterior. Then by taking into account \eqref{Kdef}, for the imaginary part of the double layer kernel $K_\kappa$ we use regular rules while, for its real part we can rely on the singular rules for the kernel $1/r.$ The constant $c$ involved in the near singularity detection phase in this experiment is chosen to be equal to $0.1$.
 In all the tests reported for this and also in the next experiment for all the three kinds of quadrature rules the QI degree is $2$ for the singular and regularized rules and it is $4$ for the regular ones. We observe that a low degree is reasonable for regularized integrals, since higher degrees of the quasi interpolating spline are profitable only if the  function it approximates is sufficiently regular. For singular integrals a low degree has been preferred to reduce the cost of the related weight computation. For all the  involved quadrature rules, in this experiment the chosen number of QI nodes (uniformly distributed in the support of each trial B-spline function) is $7,9,11,$ respectively for degrees $d = d_1=d_2 = 2,3,4,$ on both coordinate directions. Note that for any degree this corresponds to select the nodes at the knots of the trial B-spline and at their midpoints in each coordinate.

 \begin{figure}[thb]
\centering
 \subfigure[Modulus of the exact and of the numerical solution  on ${\cal C}$ for $\kappa=1.$]{
\includegraphics[trim = 0.6cm 0cm .75cm 0cm, clip = true, height=5cm]{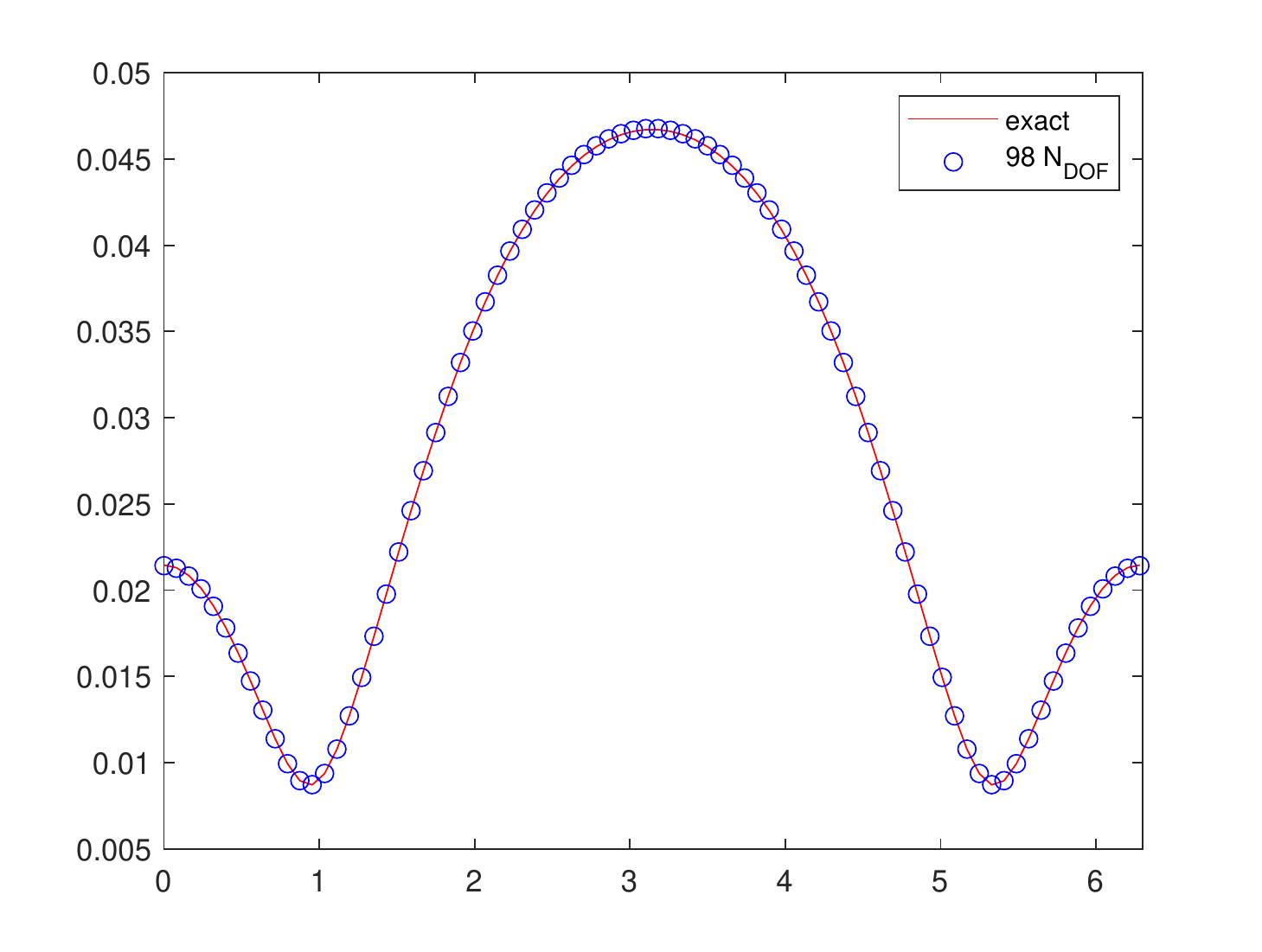}
}
\hspace{0.3cm}
\subfigure[Modulus of the exact and of the numerical solution  on ${\cal C}$ for $\kappa=3.$]{
\includegraphics[trim = 0.6cm 0cm .75cm 0cm, clip = true, height=5cm]{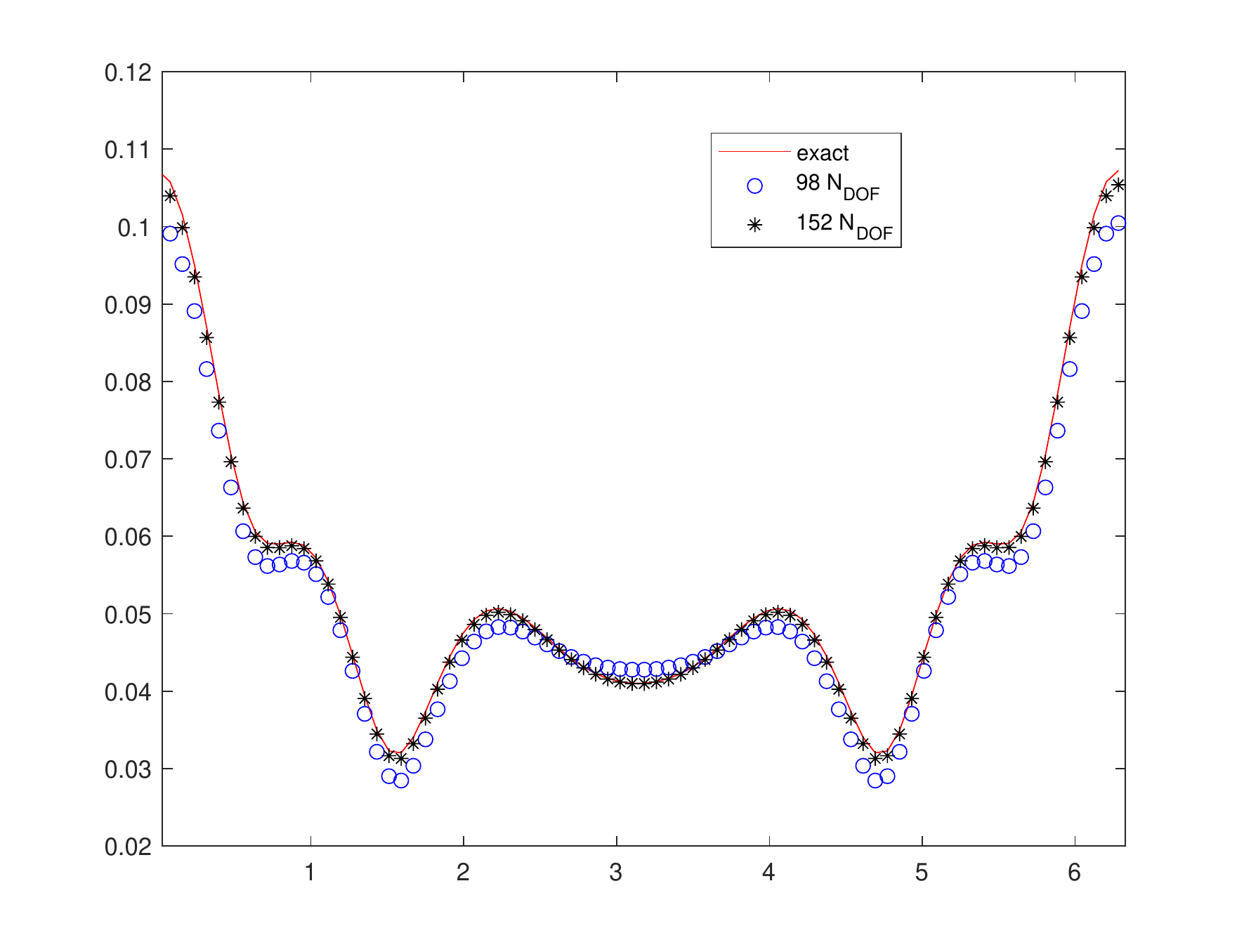}
}
\vspace{2cm}
 \subfigure[Distribution of the error $e_P$ on ${\cal C}$ for $\kappa=1.$]{
\includegraphics[trim = 0.6cm 0cm .75cm 0cm, clip = true, height=5cm]{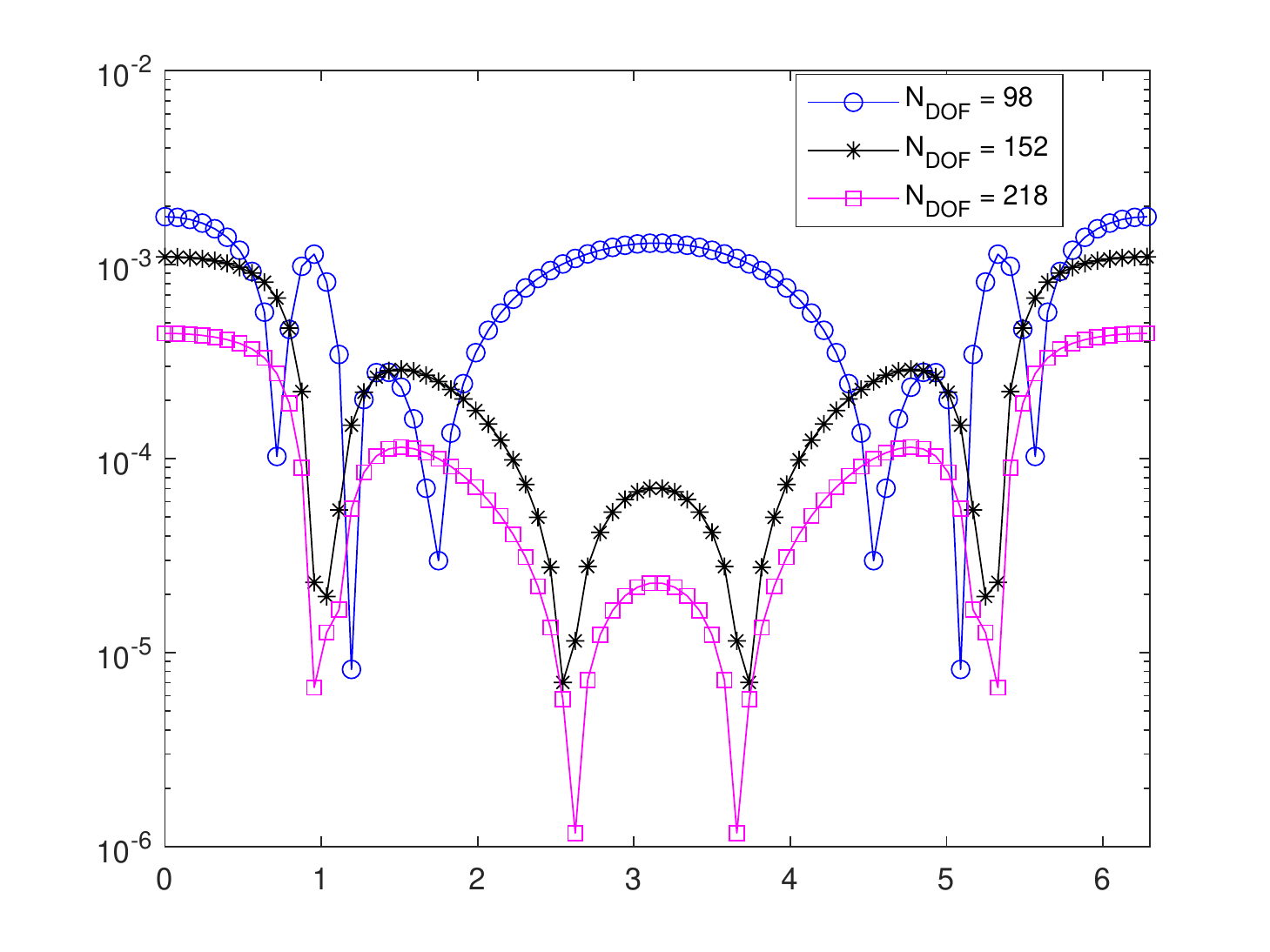}
}
\hspace{0.3cm}
 \subfigure[Distribution of the error $e_P$ on ${\cal C}$ for $\kappa=3.$]{
\includegraphics[trim = 0.6cm 0cm .75cm 0cm, clip = true, height=5cm]{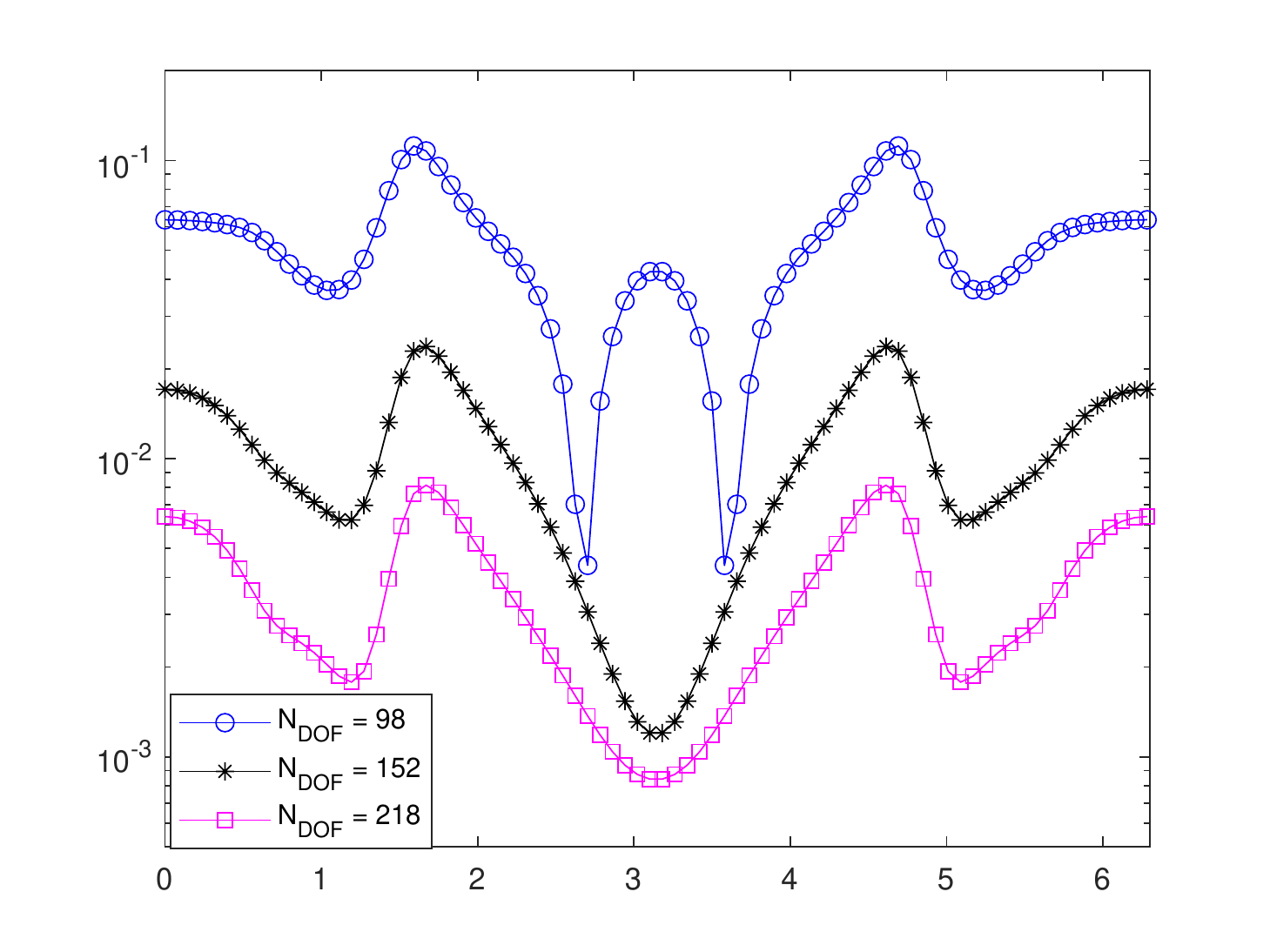}
}

\vspace{-2cm}
\caption{Rigid scattering on a sphere: results on the exterior equatorial circumference ${\cal C}$  ($d = 2, p=2,2,4$). }
\label{fig:RS_error_Simpson}
\end{figure}

In order to underline the oscillating nature of the solution of (\ref{rigidscat}) and of the related domain solution, Figure \ref{fig:SS_sol} shows the real and the imaginary part of the exact Cauchy datum $\phi = p_{tot}|\Gamma$ for $\kappa=2$, (frequency of about $100$ Hz) together with the corresponding real and imaginary parts of the total pressure in an area of the equatorial plane exterior to the sphere.  Figure \ref{fig:SS_error} shows the distribution on the sphere of the modulus of the absolute error obtained with both the choices $d=2$ and  $d=4.$  Note that the number of degrees of freedom ($N_{\rm DOF} = 218$) is the same on the left and on the right of this figure, since the number $n^2$ of elements  in each patch is chosen in order to ensure the fulfillment of such goal, that is $n=5$ for $d=2$ and $n =3$ for $d=4$. This means that when $d = 4$, the elements are larger than in the other case and so an analogous maximal error on the surface is acceptable. Note also that the error distribution on the sphere implies a slightly lower $L^2(\Gamma)$ relative error on surface for $d=4$ (it is $3.20e-03$ for $d=2$ and $2.43e-03$ for $d=4$). Our results can be compared with those, less accurate,  shown in Figure 13 of \cite{Anitescu21}, which refers to $\kappa=2$, $R=1$, $d=4$ as well, and are obtained with $N_{\rm DOF}=200$.

The off surface quality of our approach is  clarified by the plots reported in Figure \ref{fig:RS_error_Simpson}, as it is done in \cite{Simpson14}. Focusing on the approximation of the modulus of the solution, the figure shows the results obtained on the exterior equatorial circumference ${\cal C}$ with radius $10$ setting $d = 2$ (the other parameters are chosen as already described),  where the reported error is the following one
\begin{equation}
\label{eP}
e_P (\x) := \frac{ \vert \, \vert u_h(\x) \vert - \vert u^{ex}(\x) \vert \, \vert }{\vert u^{ex}(\x) \vert }.
\end{equation}
The figure outlines a good conformability between the moduli of the exact and of the numerical solution obtained with very few degrees of freedom for both $\kappa=1$ and $\kappa=3.$ The good quality of our exterior reconstructions can be further checked looking  at the  distribution on ${\cal C}$ of the error $e_P.$ Figures~\ref{fig:RS_error_Simpson}(c) and \ref{fig:RS_error_Simpson}(d) can be compared respectively with Figure~17 (where $\kappa=2$, but $R=0.5$) and Figure~20 (where $\kappa=6$, but $R=0.5$) of \cite{Simpson14}, confirming the better performance of our approach. We should remark however, as underlined in \cite{Venas20}, that the accuracy of the obtained results in \cite{Simpson14} is probably also influenced by the adopted singular parameterization of the sphere.
\begin{figure}[bht]
\centering
 \subfigure[$\kappa=1$]{
\includegraphics[trim = 0.6cm 0cm .75cm 0cm, clip = true, height=4.5cm]{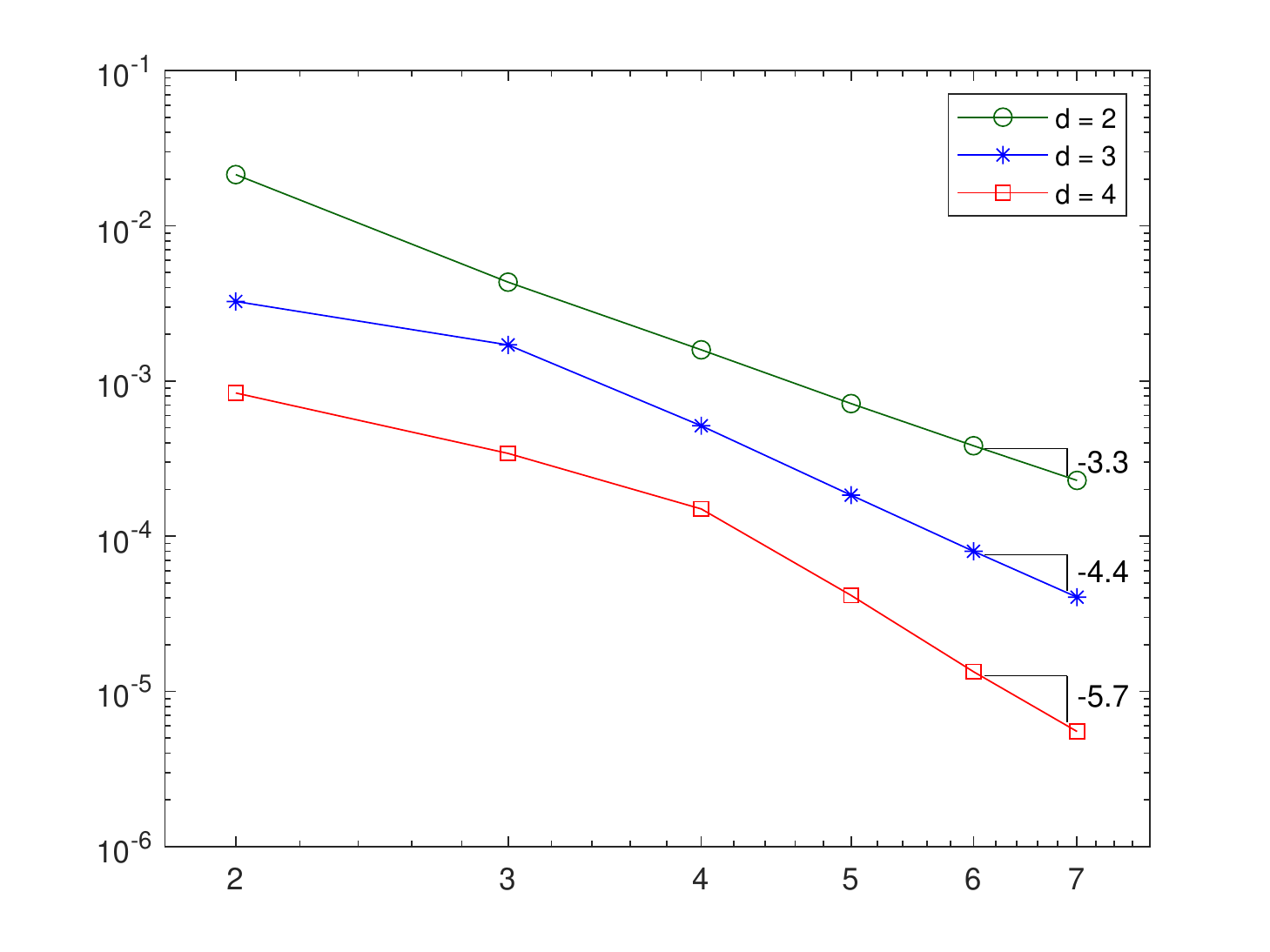}
}
\hspace{0.3cm}
\subfigure[$\kappa=2$]{
\includegraphics[trim = 0.6cm 0cm .75cm 0cm, clip = true, height=4.5cm]{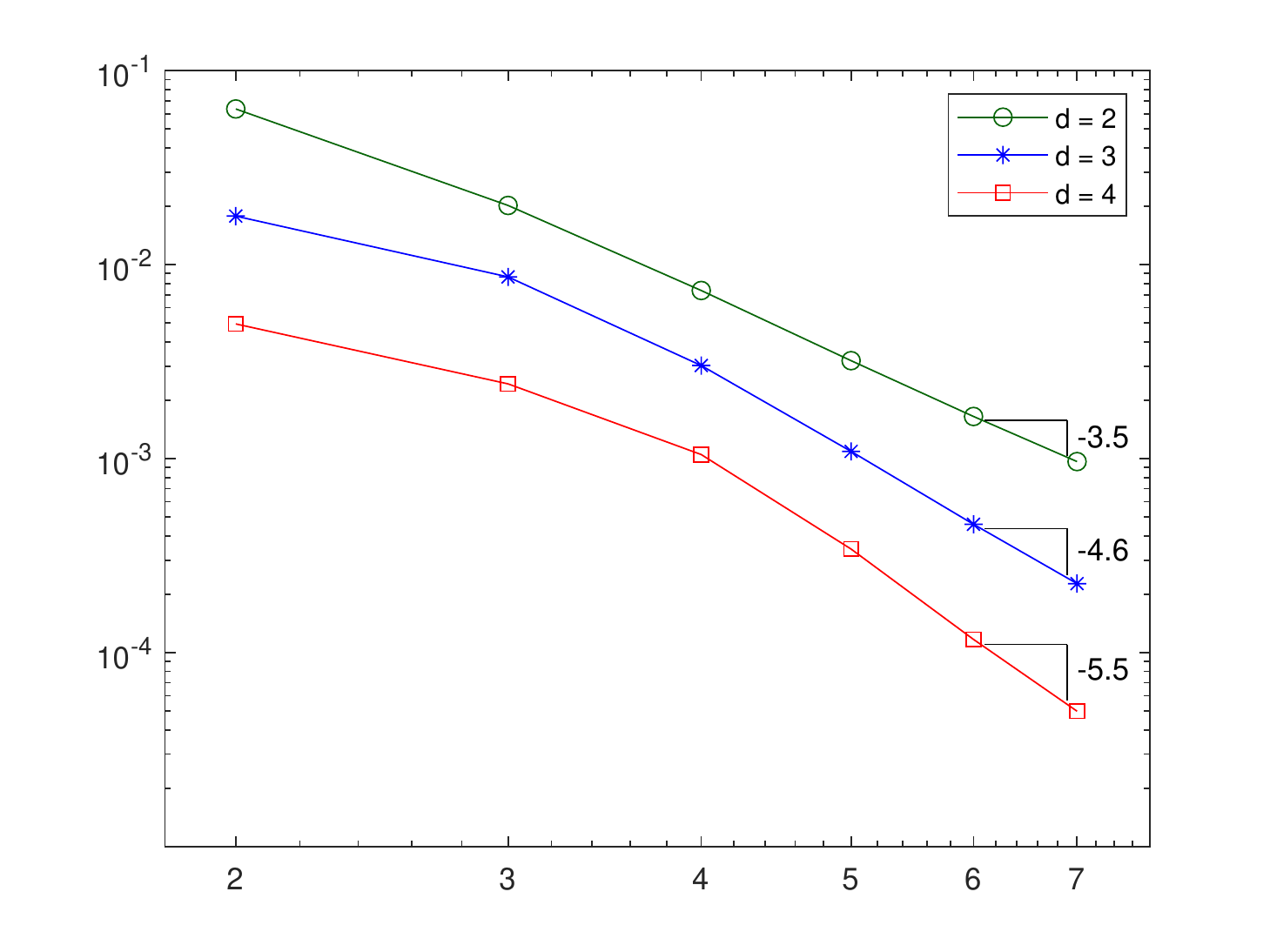}
}
\hspace{0.3cm}
 \subfigure[$\kappa=3$]{
\includegraphics[trim = 0.6cm 0cm .75cm 0cm, clip = true, height=4.5cm]{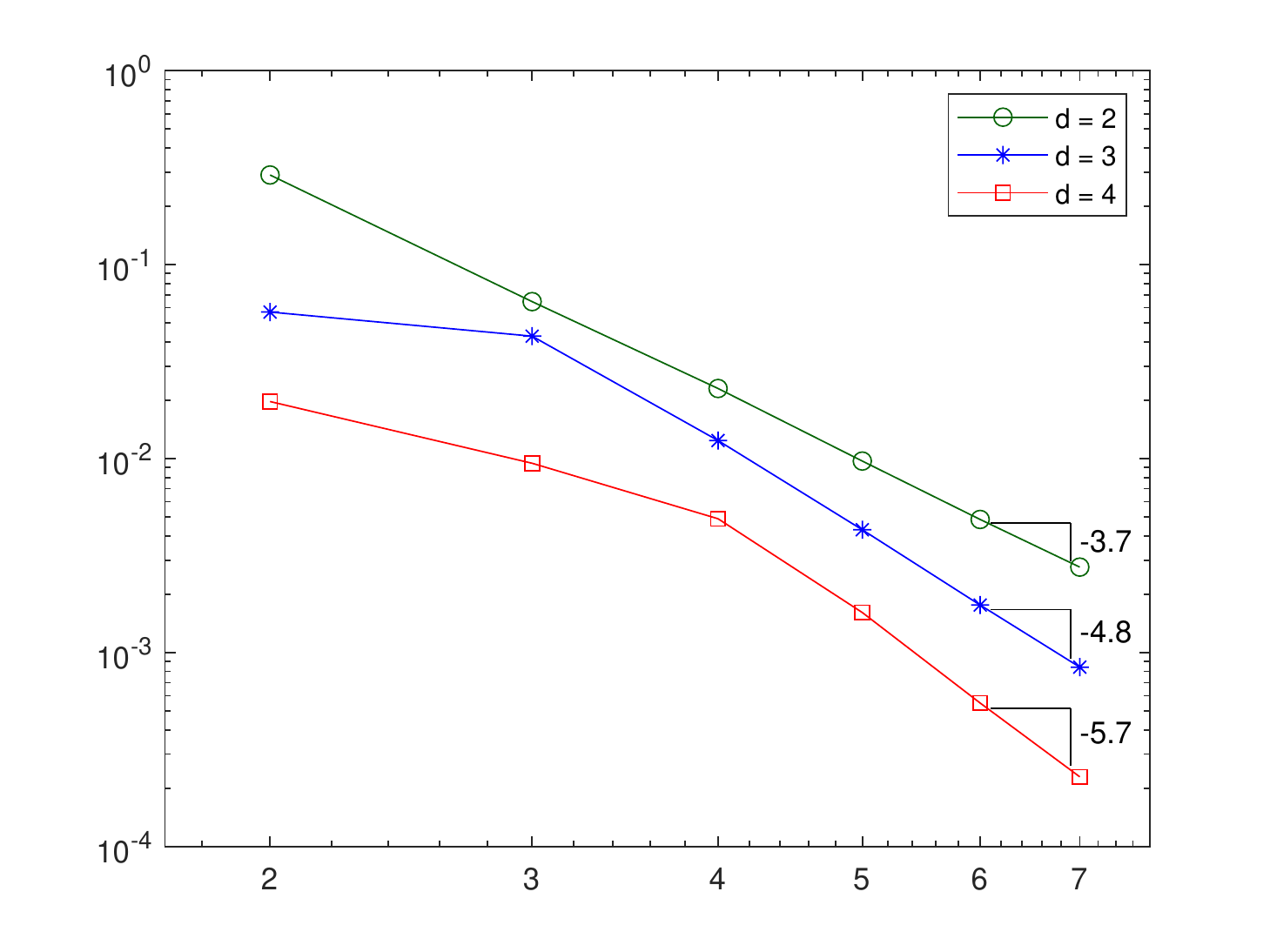}
}
\caption{Rigid scattering on a sphere: behavior of the $L^2(\Gamma)$ error, $e_{L^2}$,  versus the number $n$ of elements in each
coordinate direction for each patch.}
\label{fig:RS_error}
\end{figure}

In order to complete the analysis of the results for this experiment,  we present also Figure \ref{fig:RS_error}, which shows the dependency  from $n$  of the $L^2(\Gamma)$ relative error, $e_{L^2}$ for $d$ ranging  between $2$ and $4,$ confirming its nice convergence behavior  for all the considered values of $\kappa$. We remark that these results confirm a better performance than the approach used in \cite{Venas20}, compare the plot corresponding to $d=4$ in our Figure~\ref{fig:RS_error}(a) with the plot related to Parameterization 2 (our 6 patch parameterization of the sphere), IGABEM CCBIE (Conventional BIE, Collocation), shown in Figure~12 of \cite{Venas20}. For a precise comparison we note that $n=1,2,\ldots, 7$ in our case corresponds to $N_{\rm DOF}=98,152,218,296,386,488,602$.


 \subsection{Acoustic problem interior to a torus}

 This last example is taken from \cite{Simpson14}. The Helmholtz equation is considered with exact acoustic potential chosen to be
 \begin{align}
 \phi(\x) = \sin(\kappa x/\sqrt{3})\sin(\kappa y/\sqrt{3})\sin(\kappa z/\sqrt{3}).
 \end{align}
The Neumann datum is prescribed by computing the acoustic velocity field as $({{\partial \phi }/{\partial \n}})(\x)$ for $\x \in \Gamma$, where $\n$ denotes the outward unit normal vector to the surface.  In the experiment we consider the wavenumber $\kappa=2$ and global $C^0$ discretization space with spline bi-degree $(d,d)$ for $d=2,3$. The constant $c$ used in (\ref{eqn:thresh_dist}) for near singularity detection is set to $0.1$. Regarding the adopted quadrature scheme, we use QI degree $4$ for regular integrals and degree $2$ for singular and regularized integrands. For $d=2$ we set $m=2$ (number of terms in the singularity extraction) and the number of uniformly spaced quadrature nodes in each parametric direction on the support of a B-spline is set to $7$; hence nodes are just knots and knot midpoints of B-splines.  When $d=m=3,$ we use $9$ quadrature nodes for singular and regular integrals; for regularized integrals we use $13$ nodes. To obtain a more even knot span in both parametric directions of the physical space, for all 16 patches the starting mesh is described by the following knot vectors,
\begin{align*}
T_{1}^{(\ell)} = \begin{bmatrix}0 & 1/9 & 2/9 & 3/9 & 4/9 & 5/9 & 6/9 & 7/9 & 8/9 & 1 \end{bmatrix}, \qquad
T_{2}^{(\ell)} = \begin{bmatrix}0 & 1/3 & 2/3 & 1 \end{bmatrix}.
\end{align*}

While the imaginary part of the acoustic potential $\phi$ is zero on the entire boundary $\Gamma$,  the chosen wavenumber produces several oscillation in the real part of $\phi$, see Figure \ref{fig:TorusError}(a).
The error $e_{L^2}$ is plotted with respect to size $h$ in Figure \ref{fig:TorusError}(b); a convergence rate equal to $d+1$ can be observed. In Figure \ref{fig:TorusError}(c) and (d) the absolute error for the real part and the imaginary part of the unknown acoustic potential $\phi$ is shown, respectively.

 \begin{figure}[thb]
\centering
\subfigure[Acoustic potential.]{\includegraphics[trim = 0.6cm 1cm .75cm 0cm, clip = true, height=5.5cm]{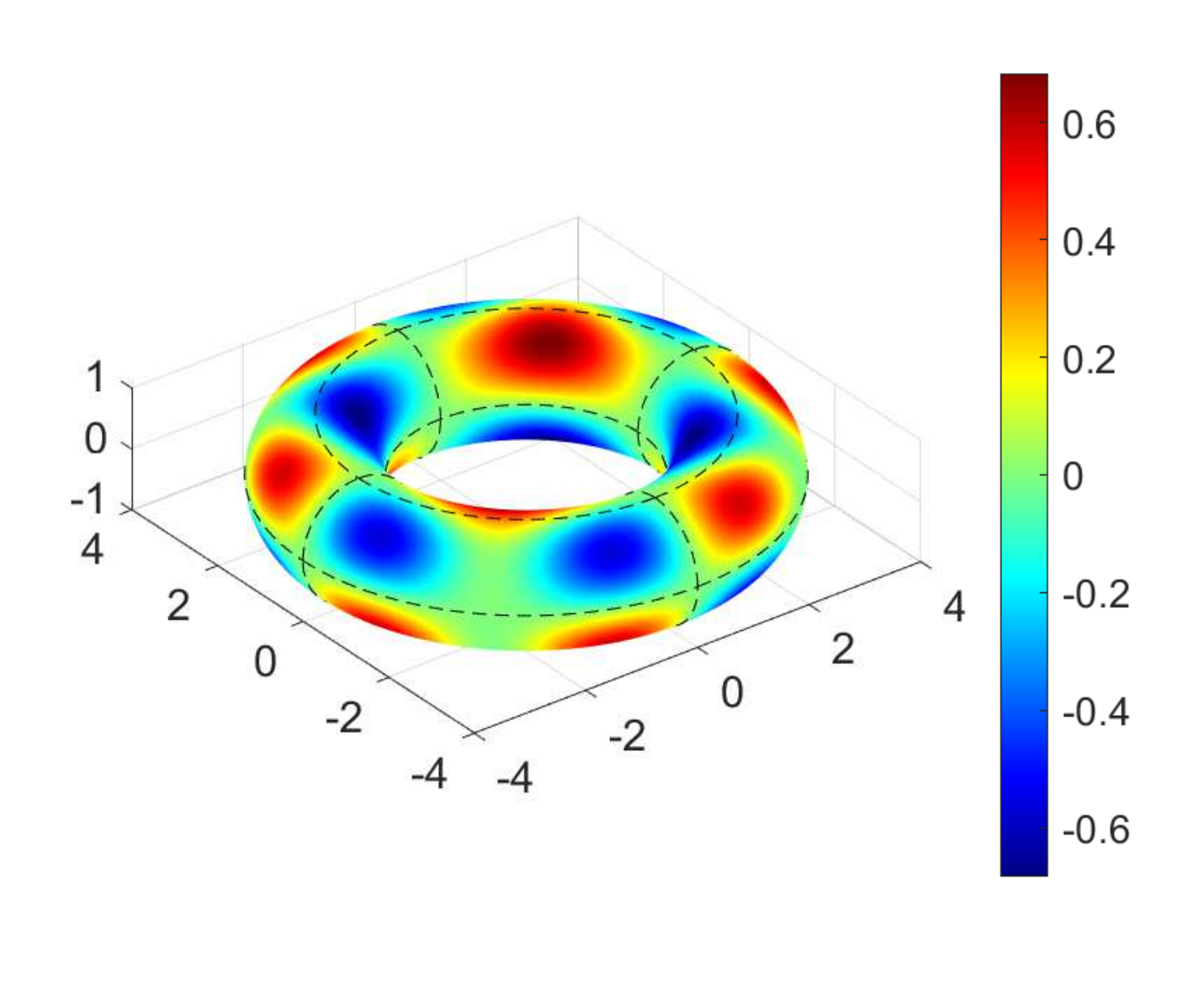}}
\qquad
\subfigure[Error convergence against mesh size $h$.] {\includegraphics[trim = 0.25cm .25cm .75cm 0.25cm, clip = true, height=5.5cm]{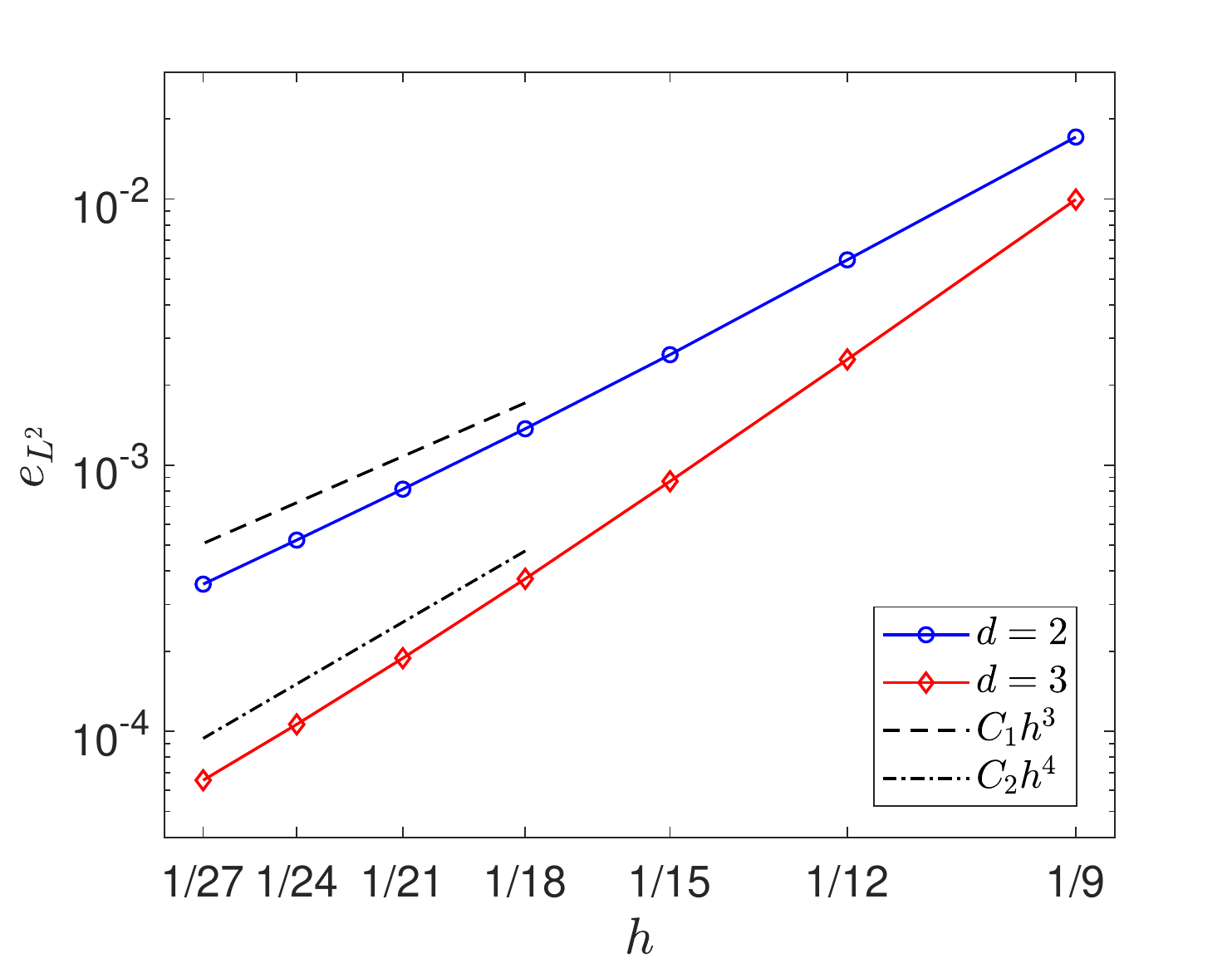}}\\
 \subfigure[On surface distribution of the real part of the absolute error for $d=3$ and $h=1/9$ ($N_{\rm DOF}=640$).]
{\includegraphics[trim = 0.6cm 1cm .5cm 0cm, clip = true, height=5.5cm]{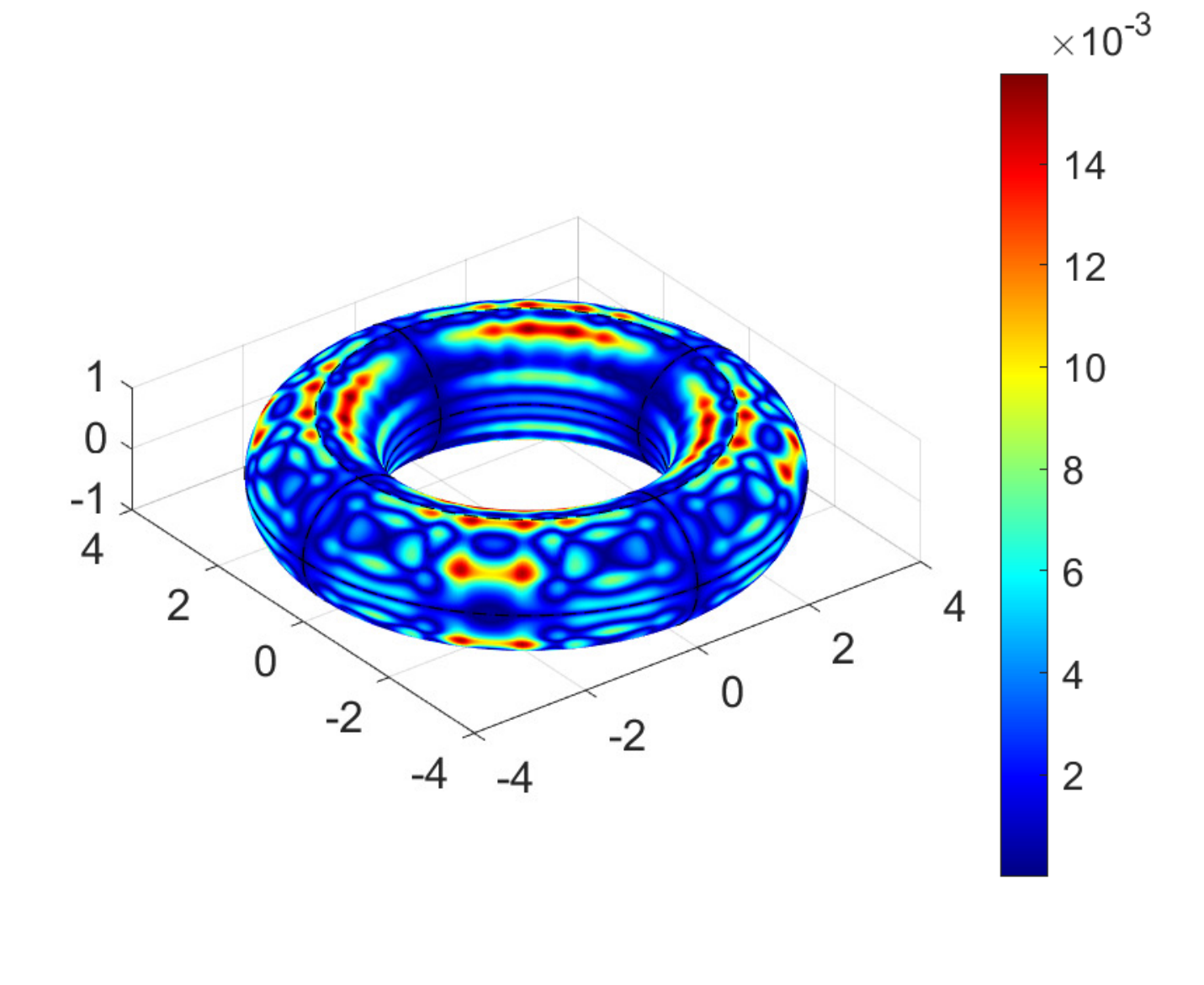}
}
\qquad
\subfigure[On surface distribution of the imaginary part of the absolute error for $d=3$ and $h=1/9$ ($N_{\rm DOF}=640$).] {\includegraphics[trim = 0.6cm 1cm .5cm 0cm, clip = true, height=5.5cm]{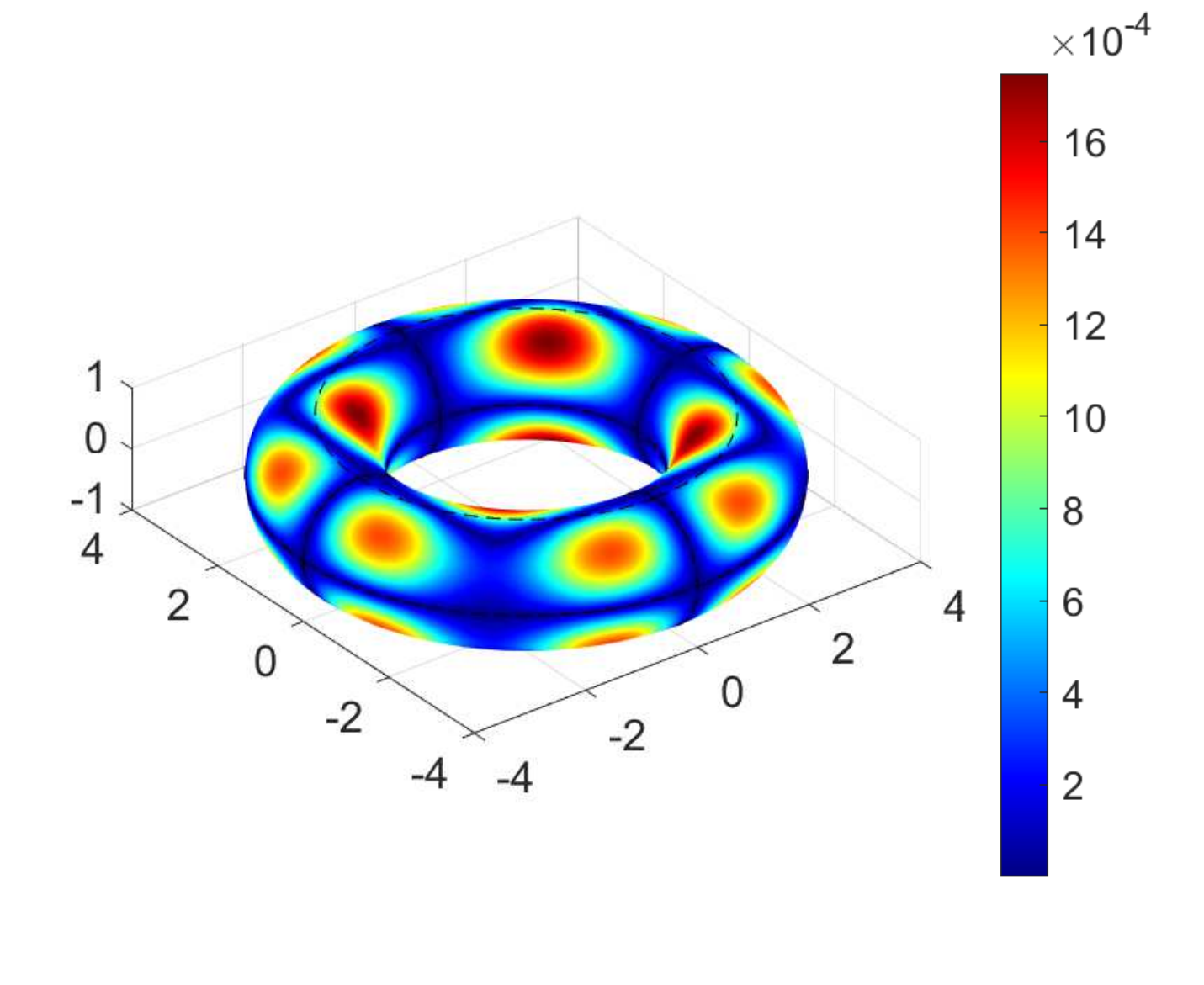}}
\caption{Acoustic problem with $\kappa=2$,  interior to a torus.}
\label{fig:TorusError}
\end{figure}

\section{Conclusion} \label{sec:Conc}

 In this paper we address 3D Helmholtz problems using isogeometric BEMs and relying on conformal multi-patch smooth  geometries and spline discretization spaces. The governing boundary integral equations are numerically approximated by a collocation scheme. The adopted quadrature rules for both regular and singular integrals are tailored for B-splines and they allow a profitable function-by-function matrix assembly. A  procedure is proposed for the automatic detection of near singular integrals. Besides studying the approximation power of the considered rules, an arrangement with a singularity extraction technique is proposed, extending it to the multi-patch setting for the first time. Numerical results confirm the effectiveness of the method by achieving sufficient accuracy of the numerical solution with a small number of uniformly distributed quadrature nodes.

A possible future work includes a generalization of the proposed simulation model to more flexible non-conformal discretization spaces and parameterizations. Further testing of the proposed methodology in both the conformal and non-conformal setting for more general types of contacts between patches could be also taken into account.

\section{Appendix} \label{sec:app}

To provide a self-contained analysis of the developed integration procedures, a more technical part related to the behavior of the considered singular kernels is given in this appendix.

Let $\Gamma$ be a regular surface admitting a $C^2$ regular parameterization on $D := [0\,,\,1]^2\,,\,\, \F: D \rightarrow \Gamma,\,\,\t=(t_1,t_2) \mapsto \F(\t)$.  
Furthermore, for any fixed $\s \in D,$ let us introduce the following notation,
$$\r_\s(\t) := \F(\t) - \F(\s)\,,\,\, r_\s(\t) := \Vert \r_\s (\t) \Vert_2 \,, \,\, \bnu(\t):= \frac{\partial \F}{\partial t_1}(\t) \times  \frac{\partial \F}{\partial t_2}(\t)\,.$$ 
\begin{prn} \label{Prop1}
For each fixed $\s \in D $ and for any unit vector $\u= (\cos \theta, \sin \theta)$ such that $\s + \rho \u$ is included in $D$ for $\rho$ positive and sufficiently small, it is
\begin{equation} \label{lim1}
\lim_{\rho \rightarrow 0^+} \frac{\rho}{r_\s(\s+\rho \u)} =  \frac{1}{\sqrt{E \cos^2 \theta +2F \sin \theta \cos \theta + G\sin^2 \theta }},    
\end{equation}
\begin{equation} \label{lim2}
\lim_{\rho \rightarrow 0^+} \frac{\r_\s(\s+\rho \u) \cdot \bnu(\s+\rho \u)}{r^2(\s+\rho \u)} = - \frac{1}{2} \frac{L \cos^2 \theta +2M \sin \theta \cos \theta + N \sin^2 \theta}{E \cos^2 \theta +2F \sin \theta \cos \theta + G \sin^2 \theta},    
\end{equation}
which are both finite limits, with $E,F,G$ and $L,M,N$ respectively denoting the coefficients of the first and second fundamental forms (with respect to the non-normalized normal $\bnu$) of $\F$ evaluated at $\s$,
\begin{align}
\label{EFG}
& E :=  \frac{\partial \F}{\partial t_1}(\s) \cdot \frac{\partial \F}{\partial t_1}(\s), && F := \frac{\partial \F}{\partial t_1}(\s) \cdot  \frac{\partial \F}{\partial t_2}(\s), && G :=   \frac{\partial \F}{\partial t_2}(\s) \cdot  \frac{\partial \F}{\partial t_2}(\s),\\
\label{LMN}
& L :=  \bnu(\s) \cdot  \frac{\partial^2 \F}{\partial t_1^2}(\s), && M := \bnu(\s) \cdot  \frac{\partial^2 \F}{\partial t_1 \partial t_2}(\s), && N :=  \bnu(\s) \cdot  \frac{\partial^2 \F}{\partial t_2^2}(\s).
\end{align} 
\end{prn}
\pf  Setting for brevity  $\F_i :=\frac{\partial \F}{\partial t_i}(\s),\,\, i=1,2$  and $\F_{ i,j} :=   \frac{\partial^2 \F}{\partial t_i \partial t_j}(\s), \, i,j= 1,2,$ and using the short notation $\r = \r(\s+\rho \u), r = r(\s+\rho \u), \bnu = \bnu(\s+\rho \u),$ we can write,
$$ \r = \rho( \F_1 \cos \theta + \F_2 \sin \theta) + \rho^2( \frac{1}{2}\F_{1,1} \cos^2 \theta + \F_{1,2}\cos \theta \sin \theta +  \frac{1}{2}\F_{2,2} \sin^2 \theta) + O(\rho^3) $$
and
$$\begin{array}{ll} \bnu &= [\F_1 + \rho \big(\F_{1,1} \cos \theta + \F_{1,2} \sin \theta \big) + O(\rho^2)] \times   [\F_2 + \rho (\F_{1,2} \cos \theta + \F_{2,2} \sin \theta) + O(\rho^2)]  = \cr
\ & = (\F_1 \times \F_2) + \rho [ (\F_1 \times \F_{1,2} - \F_2 \times \F_{1,1}) \cos \theta + (\F_1 \times \F_{2,2} - \F_2 \times \F_{1,2}) \sin \theta ] + O(\rho^2)\,. \end{array}$$
Since $ r^2  = \rho^2 (E \cos^2 \theta +2F \sin \theta \cos \theta + G\sin^2 \theta) + O(\rho^3)\,,$ we have,
$$\frac{\rho}{r} = \frac{1}{\sqrt{E \cos^2 \theta +2F \sin \theta \cos \theta + G\sin^2 \theta + O(\rho) } }\,, \quad  
\r \cdot \bnu = -\frac{1}{2}\rho^2 \big(L \cos^2 \theta +2M \sin \theta \cos \theta + N\sin^2 \theta \big) + O(\rho^3)\,.$$
In order to confirm that the two limits in (\ref{lim1}) and (\ref{lim2}) are always finite, we observe that their denominator can be respectively written as $\sqrt{\u^T M_\s \u}$ and $\u^T M_\s \u\,,$ with $M_\s$ denoting the following symmetric $2 \times 2$  positive definite matrix associated to the first fundamental form of $\F$ evaluated in $\s,$
\begin{equation} \label{matM}
M_\s := \left(\begin{array}{ll} E & F \cr F & G \end{array} \right)\,.  
\end{equation}
Now the proof is complete. \eop 
 Note that the previous result implies that both the functions  $\t \rightarrow \Vert \t - \s \Vert_2 /r_\s(\t)$ and $\t \rightarrow (\r_\s(\t) \cdot \nu(\t))/r_\s(\t)^2$ are discontinuous at $\t = \s.$ Despite this, the following corollary allows one to verify their uniform boundedness in $D$.
 \begin{crl}
 Under the assumptions of the previous proposition, both the functions $g_{1,\s}(\t) := \Vert \t - \s \Vert_2 /r_\s(\t)$ and $g_{2,\s}(\t) := (\r_\s(\t) \cdot \nu(\t))/r_\s^2(\t)$ are bounded in their domain $D \setminus \{\s\},$ uniformly with respect to $\s \in D.$
 \end{crl}
 \pf
Since the proof is analogous for $g_{1,\s}$ and $g_{2,\s},$ let us refer only to $g_{1,\s}$. 

Using proof by contradiction, we can assume that there exists a convergent sequence $\{(\s_{k_1},\t_{k_1})\}$ in $D^2$ such that $g_{1,\s_{k_1}}(\t_{k_1}) \rightarrow + \infty$. Then, considering that $\F$ is continuous and injective, it holds $\rho_{k_1} := \Vert \t_{k_1} - \s_{k_1} \Vert_2 \rightarrow 0$ and the limit of the convergent sequence is equal to $(\s,\s)$ for some $\s \in D$. Now, setting $\t_{k_1} = \s_{k_1} + \rho_{k_1}(\cos \theta_{k_1}, \sin \theta_{k_1}),$ and by considering that $\theta_{k_1}$ varies in the periodic interval $[0,2\pi]$, we can find a subsequence $\{(\s_{k_2},\t_{k_2})\}$ of $\{(\s_{k_1},\t_{k_1})\}$ such that $\{\theta_{k_2}\}$ tends to a certain angle $\theta$. Then, using arguments analogous to the ones used in the proof of the previous proposition, it can be proven that
\begin{align}\label{limit}
\lim_{k_2 \rightarrow \infty} g_{1,\s_{k_2}}(\t_{k_2}) = \frac{1}{\sqrt{E \cos^2  \theta +2 F \sin \theta \cos \theta + G\sin^2 \theta }},
\end{align}
where $E, F, G$ are defined in \eqref{EFG}.

Since the limit in \eqref{limit} exists, we obtain a contradiction, since $g_{1,\s_{k_3}}(\t_{k_3}) \rightarrow + \infty$ for every subsequence $\{(\s_{k_3},\t_{k_3})\}$ of sequence $\{(\s_{k_1},\t_{k_1})\}$. Since $\F$ is regular, we can  bound $ E,  F, G$ for all $ \s \in D$ and obtain a common upper bound for the limit  \eqref{limit} for all $\s \in D$. \eop 
 
The following last proposition  is useful to estimate an upper bound for our quadrature rules for singular integrals.
\begin{prn} \label{Prop2}
Let  $R \subset D$ be a rectangular integration domain with edges of dimensions $H_1$ and $H_2.$  Then there exists a constant $C$ such that  for each $\s \in R$ it is
\begin{equation} \label{limited}
   \int_R \vert U(\s,\t) \vert d\t \,\, \le C  \max \{H_1,H_2\} \,,
\end{equation}
where $U$ is one of the two kernels $U_{\rm SL}$ and $U_{\rm DL}$ defined in (\ref{eqn:kernelU}).
\end{prn}
\pf
First of all, setting $\rho_\s(\t) := \Vert \t - \s \Vert_2,$ we observe that we can write $U_{\rm SL}(\s,\t) =  g_{1,\s}(\t)/\rho_\s(\t)$ and
$U_{\rm DL}(\s,\t) =  g_{1,\s}(\t) g_{2,\s}(\t)/ \rho_\s(\t)\,,$ where $g_{1,\s}$ and $g_{2,\s}$ are the two functions defined in the previous corollary. Thus we can derive from the corollary itself that, for both the considered definitions of $U,$ there exists a positive constant $Q$ such that 
$$\vert U(\s,\t) \vert \le Q/\rho_\s(\t) \,, \,\, \forall (\s,\t) \in D^2,\,\, \s \ne \t\,.$$ 
Then, denoting with  ${\cal A}_\s$ the circle centered at $\s$ with radius $\sqrt{2} \max\{H_1,H_2\}$ we have that $R \subset D \cap A_\s$  which implies the following inequalities, 
$$ \int_R \vert U(\s,\t) \vert d\t \,\, \le  \int_{D \cap {\cal A}_\s} \vert U(\s,\t) \vert d\t \le Q \int_{{\cal A}_\s} \frac{1}{\rho_\s(\t)} d\t\,.$$
Thus, using polar coordinates centered at $\s$ the thesis easily descends. \eop
The constant $C$ in \eqref{limited} depends on $\F$. By prescribing common bounds for coefficients in \eqref{EFG} and \eqref{LMN} for all ``well-behaved'' functions $\F$ considered in practical applications, a common constant $C$ can be sufficiently well estimated.

\section*{Acknowledgements}
Tadej Kandu\v{c} was partially funded for developing this research by SIMAI-ACRI, which is gratefully acknowledged. The research of Antonella Falini was founded by PON Project AIM 1852414 CUP H95G18000120006 ATT1.
The Italian authors are members of the INdAM Research group GNCS. The INdAM support through GNCS and Finanziamenti Premiali SUNRISE is gratefully acknowledged.

\bibliography{BEMquad}


\end{document}